\documentclass[pdflatex]{sn-jnl}

\usepackage{graphicx}%
\usepackage{multirow}%
\usepackage{amsmath,amssymb,amsfonts}%
\usepackage{amsthm}%
\usepackage[title]{appendix}%
\usepackage{xcolor}%
\usepackage{textcomp}%
\usepackage{manyfoot}%
\usepackage{booktabs}%
\usepackage{algorithm}%
\usepackage{algorithmicx}%
\usepackage{algpseudocode}%
\usepackage{listings}%
\usepackage{tikz-cd}
\usepackage[utf8]{inputenc}
\usepackage{enumitem}

\theoremstyle{thmstyleone}%
\newtheorem{theorem}{Theorem}%
\newtheorem{proposition}[theorem]{Proposition}%
\newtheorem{lemma}[theorem]{Lemma}

\theoremstyle{thmstyletwo}%
\newtheorem{remark}{Remark}%

\theoremstyle{thmstylethree}%
\newtheorem{definition}{Definition}%

\def\endpf{\hbox{\vrule height1.5ex width.5em}}

\begin{document}

\title[Invariant differential operators on homogeneous space]{\bf Algebraic and analytic properties of invariant differential operators on a homogeneous space of complexity $1$}


\author[1]{\fnm{Hanlong} \sur{Fang}}\email{hlfang@pku.edu.cn}

\author[2]{\fnm{Xiaocheng} \sur{Li}}\email{lixiaocheng@sdu.edu.cn}

\author[3]{\fnm{Yunfeng} \sur{Zhang}}\email{zhang8y7@ucmail.uc.edu}

\affil[1]{\orgdiv{School of Mathematical Sciences}, \orgname{Peking University}, \orgaddress{ \city{Beijing}, \postcode{100871}, \country{China}}}

\affil[2]{\orgdiv{Data Science Institute}, \orgname{Shandong University}, \orgaddress{\city{Jinan}, \postcode{250100}, \state{Shandong}, \country{China}}}

\affil[3]{\orgdiv{Department of Mathematical Sciences}, \orgname{University of Cincinnati}, \orgaddress{\city{Cincinnati}, \postcode{45221}, \state{Ohio}, \country{United States}}}


\abstract{Denote by $SL_3(\mathbb R)$ the special linear group of degree 3 over the real numbers, $A$ the subgroup consisting of the diagonal matrices with positive entries. In this paper, we study the algebraic and analytic properties of the invariant differential operators on the homogeneous space $SL_3(\mathbb R)/A$. Firstly, we specify the noncommutative algebra of invariant differential operators in terms of generators and their relations. Secondly, we describe the center of this algebra and prove that all of its symmetric elements are essentially self-adjoint. Thirdly, for the first time on homogeneous spaces, we identify several essentially self-adjoint invariant differential operators which do not lie in the center of the algebra of invariant differential operators. 
}




\maketitle

\section{Introduction}

Denote by $SL_n(\mathbb R)$ the group of $n\times n$ real matrices with determinant one, $A$ the subgroup consisting of diagonal matrices with positive entries, and $SL_n(\mathbb R)/A$ the manifold of left cosets $gA$, $g\in SL_n(\mathbb R)$. In this paper, we will be concerned with the harmonic analysis on the homogeneous space  $SL_n(\mathbb R)/A$ when $n=3$.

As one of the simplest  homogeneous spaces of positive complexity, $SL_3(\mathbb R)/A$ manifests itself in various ways. It is naturally related to the space of non-degenerate triangles in the plane, which is introduced by Schubert \cite{Sc} from the perspective of enumerative geometry, and of which compactification has been attracting much attention \cite{Se,Rob,RSe}. More recently, it has also been applied in the theory of perverse sheaves and derived categories in algebraic geometry \cite{BKS}. 

In a series of works \cite{Zh1,Zh2,Zh3,ZZ1,Zh4}, R. Zhang has illustrated many interesting potential applications of spectral properties of $SL_3(\mathbb R)/A$ in homogeneous dynamics and Diophantine geometry. It is therefore desirable to investigate its more refined spectral properties, especially when contrasted with their spherical counterparts (complexity $0$), for which the Plancherel formulas are well-described by data on the boundary divisors of their compactifications  \cite{SV,DKKS}. 

Recall that the classical method for extracting spectral information of Riemannian symmetric spaces is through the study of symmetric invariant differential operators \cite{HC1,HC2}. This approach has subsequently been applied to decompose various representations \cite{St,Ros,Ba,OZ,DZ,PZ,MOZ,ZG1} and to establish connections between the corresponding eigenvalues and special functions \cite{Sh1,Sh2,Sh3,Sh4,Sh5,ZG3,ZG2}. For pseudo-Riemannian symmetric spaces, the analytic technique must be combined with the delicate intrinsic geometry \cite{FJ,OS,De,BS1,BS2} to overcome the difficulty caused by the non-ellipcity of the Laplace--Beltrami operators. 

Despite its central role in spectral theory, the algebra of invariant differential operators on $SL_3(\mathbb R)/A$ is largely unexplored. To date, it is only known to be noncommutative rather than a polynomial ring \cite{Kn}, and that the natural unitary representation of $SL_3(\mathbb R)/A$ in $L^2(SL_3(\mathbb R)/A)$ is tempered \cite{BK}. Meanwhile, among the analytic properties of a symmetric differential operator, one of the central problems is whether such an operator has a unique self-adjoint extension (see for instance \cite{Ga,Co} for the Hodge--Laplace--Beltrami operator, \cite{W} for the Dirac operator, and \cite{Ch} for certain first-order differential operators).  

This paper is devoted to the study of the algebraic and analytic properties of invariant differential operators on $SL_3(\mathbb R)/A$. More specifically, we will give an explicit presentation of the noncommutative algebra. Beyond the central elements, we will establish the essential self-adjointness for certain non-central generators of this algebra. The spectral decomposition of the invariant differential operators and its interaction with Plancherel theorems are left for future research. 
\medskip

Before describing our results in more detail, we first introduce certain notations. To treat in a unified way, let $G$ denote $SL_3(\mathbb R)$,  $\mathfrak g$ the Lie algebra of $G$, and $U(\mathfrak g)$ the universal enveloping algebra of the complexification of $\mathfrak g$. Denote by $C^{\infty}(G)$ the space of complex-valued smooth functions on $G$. Then, the infinitesimal action $R$ on $C^{\infty}(G)$ induced by the right regular representation of $G$, maps $U(\mathfrak g)$ into the algebra of algebraic differentials on $G$. More precisely, $R$ acts on $u=\sum X_{1}X_2\cdots X_k\in U(\mathfrak g)$ by
\begin{equation*}
\left(R(u) f\right)(g): =\left(R_{u} f\right)(g):=\sum \left.\frac{\partial}{\partial t_1}\right|_{t_1=0}\cdots\left.\frac{\partial}{\partial t_k}\right|_{t_k=0}f(g\exp \left(t_1X_1\right)\cdots\exp \left(t_kX_k\right)).
\end{equation*}
Here $X_{1}, X_2, \cdots,X_k\in\mathfrak g$, and $f \in C^{\infty}(G)$; the exponential map is given by $\exp(X):=\gamma (1)$, where
$\gamma\colon \mathbb R \to G$ is the one-parameter subgroup of $G$ whose tangent vector at the identity is equal to $X$. It is easy to verify that $R_u$, $u\in U(\mathfrak g)$, is a left $G$-invariant differential operator on $G$. Denote by $\mathbb D(G)$ the algebra of the left $G$-invariant differential operators on $G$.

For a closed subgroup $H\subset G$, denote by $\mathbb D(G/H)$ the algebra of $G$-invariant differential operators on the homogeneous space $G/H$. Denote by $\pi:G\rightarrow G/H$ the natural projection and $\mathfrak h$ the Lie algebra of $H$. Define
\begin{equation*}
\mathbb D^H(G):=\{D\in\mathbb D(G)\,|\, D(f\circ R_h)\circ R_h^{-1}=Df,\,\forall h\in H\,\,{\rm and}\,\, f\in C^{\infty}(G)\},
\end{equation*}
where $R_h:g\mapsto gh$ is the right translation of $G$ for $h\in H$. Assuming $G$ and $H$ are reductive, we have the standard isomorphism (Theorem 4.6 in Chapter 2 of \cite{He})
\begin{equation}\label{hgh}
\mathbb D^H(G)/\left(\mathbb D^H(G)\cap \mathbb D(G)\mathfrak h\right)\cong\mathbb D(G/H).
\end{equation}
It is induced by the map 
$\mu:\mathbb{D}^H(G)\to \mathbb{D}(G/H)$, 
such that 
for each $D\in \mathbb D^H(G)$, $\mu(D)$ is the element of $\mathbb D(G/H)$ such that 
\begin{equation}\label{mu}
(\mu(D)f)\circ\pi=D(f\circ\pi)\, \,\, {\rm for\,\,all\,\,smooth\,\,functions\,\,}f\,\,{\rm on\,\,}G/H.
\end{equation}

Let $S(\mathfrak g)$ be the symmetric algebra over $\mathfrak g$. Then for a basis $\{X_1,\cdots,X_n\}$ of $\mathfrak g$, $S(\mathfrak g)$ can be identified with the algebra of polynomials
\begin{equation*}
	\sum\nolimits_{(k_1,\cdots,k_n)\in\mathbb N^{n}}a_{k_1\cdots k_n}X_1^{k_1}\cdots X_n^{k_n},\,\,a_{k_1\cdots k_n}\in\mathbb C.
\end{equation*}
We have the following  symmetrizer map
$\lambda:S(\mathfrak{g})\to \mathbb{D}(G)$. 
\begin{theorem}[Theorem 4.3 in Chapter 2 of \cite{He}]\label{HCiso} There is a unique linear bijection $\lambda$ from $S(\mathfrak g)$ to $\mathbb D(G)$ such that $\lambda(X^m)=R(X^m)$ for $X\in\mathfrak g$ and $m\in\mathbb Z^+$. More precisely,  
\begin{equation*}
	(\lambda(P)f)(g):=\left.P\left(\frac{\partial}{\partial t_1},\cdots,\frac{\partial}{\partial t_n}\right)f\left(g\exp(t_1X_1+\cdots +t_nX_n)\right)\right|_{t_1=\cdots =t_n=0},
\end{equation*}
for $P\in S(\mathfrak g)$ and $f\in C^{\infty}(G)$. In particular (see Page 282 of \cite{He}), 
\begin{equation*}
\lambda(Y_1\cdots Y_k) =\frac{1}{k!}\sum_{\sigma\in S_k}R\left(Y_{\sigma(1)}\cdots Y_{\sigma(k)}\right),
\end{equation*}
where $Y_1,\cdots, Y_k\in\mathfrak g$, and $S_k$ is the symmetric group of degree $n$. 
\end{theorem}

Denote by $E_{ij}$ the $3\times 3$ matrix unit with a $1$ in the $i^{\rm th}$ row and $j^{\rm th}$ column. For distinct $i,j,k\in\{1,2,3\}$, define differential operators on $SL_3(\mathbb R)/A$
\begin{equation}\label{d12d123}
\begin{split}
&D_{ij}:=(\mu\circ\lambda)\left(E_{ij}E_{ji}\right),   \,\,(i,j)=(1,2),(1,3),(2,3),\\
&D_{ijk}:=(\mu\circ\lambda)\left(E_{ij}E_{jk}E_{ki}\right),\,\,(i,j,k)=(1,2,3),(2,1,3).    
\end{split}
\end{equation}
Note that $D_{ij}=D_{ji}$ and $D_{ijk}=D_{jki}=D_{kij}$. We prove that
\begin{theorem}\label{main1}
$\mathbb D\left(SL_3(\mathbb R)/A\right)$ is the noncommutative associative algebra generated over $\mathbb C$ by $\{D_{12}$, $D_{13},D_{23},D_{123},D_{213}\}$ with relations
\begin{equation}
\left\{\begin{aligned}
&[D_{123},D_{213}]=0,\\
&[D_{ij},D_{ik}]=D_{ijk}-D_{ikj},\,\,\,\,\,\,i,j,k\in\{1,2,3\}{\rm\,\,are\,\,distinct},\\			&[D_{ijk},D_{ij}]=D_{jk}D_{ij}-D_{ij}D_{ik},\,\,\,\,\,\,i,j,k\in\{1,2,3\}{\rm\,\,are\,\,distinct}, \\
&2\left(D_{123}D_{213}+D_{213}D_{123}-D_{12}D_{23}D_{31}-D_{13}D_{32}D_{21}\right)=\left(D_{23}-D_{13}-D_{12}\right)^2.\\
\end{aligned}\right.
\end{equation}	
The center of $\mathbb D\left(SL_3(\mathbb R)/A\right)$ is a polynomial ring in $D_{123}+D_{213}$ and $D_{12}+D_{23}+D_{13}$.
\end{theorem}

For general homogeneous spaces, a classical result shows that the symmetric elements in the image of the center of the universal enveloping algebra are essentially self-adjoint \cite{Seg,NS}). Then, 
\begin{proposition}\label{easy}
Every symmetric differential operator in the center of $\mathbb D\left(SL_3(\mathbb R)/A\right)$ is essentially self-adjoint.    
\end{proposition}

For elements not lying in the image of the center of the universal enveloping algebra, van den Ban \cite{Ba} established essential self-adjointness of all symmetric invariant differential operators on any semisimple symmetric space. Beyond that, to the best of our knowledge, there is no general theory ensuring the essential self-adjointness in the pseudo-Riemannian setting, even for the Laplacian operators (see \cite{KK}).

The major part of the paper is  devoted to
\begin{theorem}\label{self}
The differential operators $D_{12},D_{13},D_{23}$ on $SL_3(\mathbb R)/A$ are essentially self-adjoint.
\end{theorem}

We now briefly describe the basic ideas for the proofs. To determine the algebraic structure of $\mathbb D(SL_3(\mathbb R)/A)$, the Huang-Yin normal form theory \cite{HY} plays an essential role, which solves a large system of linear equations in an inductive way. To study the essential self-adjointness of symmetric operators, we modify the scheme of \cite{Ba}. Recall that the elegant proof in \cite{Ba} is to decompose the invariant differential operator into a bounded sum of left derivatives so that the wild growth of the coefficients can be treated as bounded ones. Unfortunately, in this non-spherical case, the left derivatives are too degenerate to span the whole space of invariant differentials in a mild way. We make the observation that by choosing the cutoff functions and the mollifiers in a compatible way instead of isolating the G\r arding type space as the operator core, one may gain extra control of the coefficients. In fact, the chosen cutoff functions are annihilated by the wildest terms and contribute the desired decays thereafter.

The organization of the paper is as follows. \S \ref{stru} is devoted to the algebraic structure of $\mathbb D(SL_3(\mathbb R)/A)$. In \S \ref{central}, we prove Proposition \ref{easy}. After introducing a coordinate system for $SL_3(\mathbb R)/A$ in \S \ref{coordinate}, we represent $D_{12}$ as left derivatives in \S \ref{prop}. In \S \ref{extra}, we establish the density of $C^{\infty}_c\left(SL_3(\mathbb R)/A\right)$ in ${\rm Dom}(D_{12})$ in the graph norm, and, as a consequence, prove Theorem \ref{self}. The explicit formulas for the generators of the left derivatives and of the left invariant differentials  are given in Appendices \ref{LEFTD} and \ref{lefti}, respectively.

\section*{Acknowledgement}  The authors appreciate greatly Professors N.~Li, W.~Li, Y.~Xu and J.~Yu for many helpful discussion.  This research is supported by National Key R\& D Program of China (No.~2022YFA1006700). The first author is partially supported by NSFC
grant (No.~12201012).

\section{Structure of the Algebra of the Invariant Differentials}\label{stru}

\subsection{Generators and Relations}
Denote by $E_{ij}$ the $n\times n$ matrix unit with a $1$ in the $i^{\rm th}$ row and $j^{\rm th}$ column. Define \begin{equation}\label{slba}
\begin{split}
&X_{ij}:=E_{ij},\,\,1\leq i\neq j\leq n,\,\,\,\,\,X_{ll}:=E_{ll}-E_{nn},\,\, 1\leq l\leq n-1,
\end{split}    
\end{equation}
which constitute a basis of $\mathfrak {sl}_n(\mathbb R)$.


\begin{lemma}\label{generators}
The algebra	$\mathbb D(SL_n(\mathbb R)/A)$ is generated by 
\begin{equation*}
\left\{(\mu\circ\lambda)\left(E_{i_1i_2}E_{i_2i_3}\cdots E_{i_{k-1}i_k} E_{i_ki_1}\right)\left|\,\begin{matrix}
		2\leq k\leq n,\,1\leq i_1,i_2,\cdots,i_k\leq n\\
		i_1,i_2,\cdots,i_k\,\,{\rm are\,\, distinct}
	\end{matrix} \right.  \right\}.
\end{equation*}
\end{lemma}
{\bf\noindent Proof of Lemma \ref{generators}.} 
By \eqref{hgh}, the invariant differential operators on $SL_n(\mathbb R)/A$ are induced from the left $SL_n(\mathbb R)$ and right $A$ invariant differential operators on $SL_n(\mathbb R)$. Notice that $S(\mathfrak {sl}_n(\mathbb R))^A\cong \mathbb D^A(SL_n(\mathbb R))$. Then it suffices to prove that $S(\mathfrak {sl}_n(\mathbb R))^A$ is generated by 
\begin{equation*} 
\left\{X_{i_1i_2}X_{i_2i_3}\cdots X_{i_{k-1}i_k} X_{i_ki_1}\left|\,\begin{matrix}
		2\leq k\leq n,\,1\leq i_1,i_2,\cdots,i_k\leq n\\
		i_1,i_2,\cdots,i_k\,\,{\rm are\,\, distinct}
	\end{matrix} \right.  \right\},
\end{equation*}
which follows easily from the definition of the ${\rm Ad}(A)$ action on $S(\mathfrak {sl}_n(\mathbb R))$. 

The proof of Lemma \ref{generators} is complete.\,\,\,$\endpf$
\medskip


By Lemma \ref{generators} and Theorem \ref{HCiso}, the algebra of the $SL_3(\mathbb R)$-invariant differential operators on $SL_3(\mathbb R)/A$ are generated by (\ref{d12d123}).
Then, we have
\begin{proposition}\label{comm1} We have the following commutation relations. 
\begin{equation*}
[D_{12},D_{13}]=D_{123}-D_{132},\,
[D_{21},D_{23}]=D_{213}-D_{231},\,
[D_{31},D_{32}]=D_{312}-D_{321}. 
\end{equation*}
\end{proposition}
{\bf\noindent Proof of Proposition \ref{comm1}.} 
It suffices to prove the first equality. Computation yields that
\begin{equation*}
\begin{split}
&D_{123}-D_{213}=\mu\left(R\left(E_{12}E_{23}E_{31}\right)\right)-\mu\left(R\left(E_{13}E_{32}E_{21}\right)\right).\\
\end{split}    
\end{equation*}

In terms of universal enveloping algebra, we have
\begin{equation*}
\begin{split}
&E_{12}E_{21}E_{13}E_{31}=E_{13}E_{31}E_{12}E_{21}+E_{12}E_{23}E_{31}-E_{13}E_{32}E_{21},\\
&E_{21}E_{12}E_{13}E_{31}=E_{13}E_{31}E_{21}E_{12}+E_{23}E_{12}E_{31}-E_{13}E_{21}E_{32},\\
&E_{12}E_{21}E_{31}E_{13}=E_{31}E_{13}E_{12}E_{21}+E_{12}E_{31}E_{23}-E_{32}E_{13}E_{21},\\
&E_{21}E_{12}E_{31}E_{13}=E_{31}E_{13}E_{21}E_{12}+E_{31}E_{23}E_{12}-E_{21}E_{32}E_{13}.\\
\end{split}    
\end{equation*}
Also,
\begin{equation*}
\begin{split}
&E_{12}E_{31}E_{23}=E_{12}E_{23}E_{31}-E_{12}E_{21},\,\,\,E_{23}E_{12}E_{31}=E_{12}E_{23}E_{31}-E_{13}E_{31},\\
&E_{23}E_{31}E_{12}=E_{12}E_{23}E_{31}+E_{23}E_{32}-E_{13}E_{31},\,\,\,E_{31}E_{12}E_{23}=E_{12}E_{23}E_{31}+E_{32}E_{23}-E_{12}E_{21},\\
&E_{31}E_{23}E_{12}=E_{12}E_{23}E_{31}-E_{12}E_{21}+E_{32}E_{23}-E_{31}E_{13},\\
\end{split}    
\end{equation*}
and
\begin{equation*}
\begin{split}
&E_{32}E_{13}E_{21}=E_{13}E_{32}E_{21}-E_{12}E_{21},\,\,\,E_{13}E_{21}E_{32}=E_{13}E_{32}E_{21}-E_{13}E_{31},\\
&E_{21}E_{13}E_{32}=E_{13}E_{32}E_{21}+E_{23}E_{32}-E_{13}E_{31},\,\,\,E_{32}E_{21}E_{13}=E_{13}E_{32}E_{21}+E_{32}E_{23}-E_{12}E_{21},\\
&E_{21}E_{32}E_{13}=E_{13}E_{32}E_{21}-E_{12}E_{21}+E_{32}E_{23}-E_{31}E_{13}.\\
\end{split}    
\end{equation*}
Hence, 
\begin{equation*}
\begin{split}
 &D_{12}D_{13}-D_{13}D_{12}=\frac{1}{4}\mu\left(\left(R\left(E_{12}E_{21}\right)+R\left(E_{21}E_{12}\right)\right)\left(R\left(E_{13}E_{31}\right)+R\left(E_{31}E_{13}\right)\right)\right)\\ 
&\,\,\,\,\,\,\,\,\,-\frac{1}{4}\mu\left(\left(R\left(E_{13}E_{31}\right)+R\left(E_{31}E_{13}\right)\right)\left(R\left(E_{12}E_{21}\right)+R\left(E_{21}E_{12}\right)\right)\right)\\ &=\frac{1}{4}\mu\left(R\left(E_{12}E_{21}E_{13}E_{31}+E_{21}E_{12}E_{13}E_{31}+E_{12}E_{21}E_{31}E_{13}+E_{21}E_{12}E_{31}E_{13}\right)\right)\\ 
&\,\,\,\,\,\,\,\,\,-\frac{1}{4}\mu\left(R\left(E_{13}E_{31}E_{12}E_{21}+E_{13}E_{31}E_{21}E_{12}+E_{31}E_{13}E_{12}E_{21}+E_{31}E_{13}E_{21}E_{12}\right)\right)\\ 
&=\frac{1}{4}\mu\left(R\left(E_{12}E_{21}E_{13}E_{31}-E_{13}E_{31}E_{12}E_{21}\right)\right)+\frac{1}{4}\mu\left(R\left(E_{21}E_{12}E_{13}E_{31}-E_{13}E_{31}E_{21}E_{12}\right)\right)\\ 
&\,\,\,\,\,+\frac{1}{4}\mu\left(R\left(E_{12}E_{21}E_{31}E_{13}-E_{31}E_{13}E_{12}E_{21}\right)\right)+\frac{1}{4}\mu\left(R\left(E_{21}E_{12}E_{31}E_{13}-E_{31}E_{13}E_{21}E_{12}\right)\right)\\
&=\frac{1}{4}\mu\left(R\left(E_{12}E_{23}E_{31}+E_{23}E_{12}E_{31}+E_{12}E_{31}E_{23}+E_{31}E_{23}E_{12}\right)\right)\\ 
&\,\,\,\,\,\,\,\,\,\,\,\,\,\,\,-\frac{1}{4}\mu\left(R\left(E_{13}E_{32}E_{21}+E_{13}E_{21}E_{32}+E_{32}E_{13}E_{21}+E_{21}E_{32}E_{13}\right)\right)\\
&=\mu\left(R\left(E_{12}E_{23}E_{31}\right)\right)-\mu\left(R\left(E_{13}E_{32}E_{21}\right)\right).\\
\end{split}    
\end{equation*}

We complete the proof of Proposition \ref{comm1}. \,\,\,$\endpf$
\medskip

Similarly, we have
\begin{proposition}\label{comm2} We have the following commutation relations. 
\begin{equation*}
\begin{split}			[D_{123},D_{12}]=-[D_{213},D_{12}]&=D_{23}D_{12}-D_{12}D_{13},\\
[D_{312},D_{31}]=-[D_{132},D_{31}]&=D_{13}D_{12}-D_{23}D_{13},\\
[D_{231},D_{23}]=-[D_{321},D_{23}]&=D_{23}D_{31}-D_{12}D_{23}.
\end{split}  
\end{equation*}
\end{proposition}

\begin{proposition}\label{comm4} We have the following equality.
	\begin{equation}\label{2332}
\begin{split}
&D_{123}D_{213}+D_{213}D_{123}-D_{12}D_{23}D_{31}-D_{13}D_{32}D_{21}=\frac{1}{2}\left(D_{23}-D_{13}-D_{12}\right)^2.\\
\end{split}
\end{equation}
\end{proposition}

\begin{proposition}\label{comm3} We have the following commutation relation. 
\begin{equation}\label{123213}
\begin{split}
[D_{123},D_{213}]&=0.\\
\end{split}
\end{equation}
\end{proposition}
{\bf\noindent Proof of Proposition  \ref{comm2}, \ref{comm4}, \ref{comm3}.} The proof is the same as that of Proposition \ref{comm1}. We omit it here for simplicity. \,\,\,$\endpf$

\begin{remark}\label{cfree}
Let $D=\sum_{i}X_{i,1}X_{i,2}\cdots X_{i,k_i}$ be an element of $U(\mathfrak{g})$, understood as a left-invariant differential operator on $G$. Then the formal adjoint $D^*$ is given by 
\begin{equation}
D^*=\sum_{i}(-1)^{k_i}X_{i,k_i}X_{i,k_i-1}\cdots X_{i,1}.  
\end{equation}
One can thus easily verify that $D_{12},D_{13},D_{23},\sqrt{-1}\cdot D_{123},\sqrt{-1}\cdot D_{213}$ are formally self-adjoint. This is another  reason for taking $D_{ij}$, $D_{ijk}$ as the generators.
\end{remark}

\subsection{The center of the algebra of the invariant differentials}

\begin{proposition}\label{basis}A basis of the linear space $\mathbb D(SL_3(\mathbb R)/A)$ over $\mathbb C$ is given by
\begin{equation}\label{47}
\begin{aligned}
&D_{12}^kD_{23}^jD_{123}^lD_{213}^m,\,\, D_{12}^kD_{31}^iD_{123}^lD_{213}^m,\,\, D_{23}^jD_{31}^iD_{123}^lD_{213}^m,\,\,\,\, i,j,k\in\mathbb Z^+,\,\,l,m\in\mathbb Z{^{\geq 0}},\\
&D_{12}^kD_{123}^lD_{213}^m,\,\,D_{23}^jD_{123}^lD_{213}^m,\,\,D_{31}^iD_{123}^lD_{213}^m,\,\,\,\, i,j,k\in\mathbb Z^+,\,\,l,m\in\mathbb Z{^{\geq 0}},\\
&D_{123}^lD_{213}^m,\,\,l,m\in\mathbb Z{^{\geq 0}}.
\end{aligned}
\end{equation}
\end{proposition}
{\bf\noindent Proof of Proposition \ref{basis}.} According to Lemma \ref{generators},  $\mathbb D(SL_3(\mathbb R)/A)$ is a linear span of $D_{12}^kD_{23}^jD_{31}^iD_{123}^lD_{213}^m$. By (\ref{2332})  along with the commutation relations in Proposition \ref{comm1}, \ref{comm2}, and \ref{comm3}, we can replace $D_{12}D_{23}D_{31}$ by $D_{123}D_{213}$ with certain lower order terms. Hence, $D_{12}^kD_{23}^jD_{31}^iD_{123}^lD_{213}^m$, $k,j,i\geq 1$, is generated by the elements in (\ref{47}).

Next, we shall show that the elements in  (\ref{47}) are linearly independent. Suppose 
\begin{equation}\label{rel}
	\begin{split}
		0=&\sum_{\substack{l,m\in\mathbb Z{^{\geq 0}}}}r_{000lm}D_{123}^lD_{213}^m+\sum_{\substack{k\in\mathbb Z^+\\l,m\in\mathbb Z{^{\geq 0}}}}r_{k00lm}D_{12}^kD_{123}^lD_{213}^m+\sum_{\substack{j\in\mathbb Z^+\\l,m\in\mathbb Z{^{\geq 0}}}}r_{0j0lm}D_{23}^jD_{123}^lD_{213}^m\\
		&+\sum_{\substack{i\in\mathbb Z^+\\l,m\in\mathbb Z{^{\geq 0}}}}r_{00ilm}D_{31}^iD_{123}^lD_{213}^m+	
		\sum_{\substack{k,j\in\mathbb Z^+\\l,m\in\mathbb Z{^{\geq 0}}}}r_{kj0lm} D_{12}^kD_{23}^jD_{123}^lD_{213}^m\\
		&+\sum_{\substack{k,i\in\mathbb Z^+\\l,m\in\mathbb Z{^{\geq 0}}}}r_{k0ilm} D_{12}^kD_{31}^iD_{123}^lD_{213}^m+\sum_{\substack{j,i\in\mathbb Z^+\\l,m\in\mathbb Z{^{\geq 0}}}}r_{0jilm}D_{23}^jD_{31}^iD_{123}^lD_{213}^m,\\
	\end{split}
\end{equation}
where $r_{kjilm}\neq 0$ for only finitely many nonnegative integers $k,j,i,l,m$.

Let $\mathfrak p$ be the subspace of $\mathfrak {sl}_3(\mathbb R)$ spanned by $\{X_{ij}\}_{1\leq i,j\leq 3,\,i\neq j}$. Equip $\mathfrak p$ with natural coordinates $x=(x_{12},x_{13},x_{23},x_{21},x_{31},x_{32})$ such that  $\eta=\sum_{1\leq i,j\leq n,\,i\neq j}x_{ij} X_{ij}$ for each $\eta\in\mathfrak p$. Denote by $j:\mathfrak p\hookrightarrow \mathfrak g$ the natural injection. 
Define a mapping $\mathcal {P}:\mathfrak p\rightarrow SL_3(\mathbb R)/A$ by $\mathcal P:=\pi\circ\exp\circ j$.  Then, the following diagram commutes.
\begin{equation*}
	\begin{tikzcd} &\mathfrak g  \arrow{r}{\exp}&[4em]SL_3(\mathbb R) \arrow{d}{\pi} \\ &\mathfrak p\arrow[hookrightarrow]{u}{j} \arrow{r}{\mathcal P}&SL_3(\mathbb R)/A\\ 
	\end{tikzcd}\,\,.\vspace{-20pt} 
\end{equation*}

{\bf\noindent Claim.}
$\mathcal P$ is a local homeomorphism near $0\in\mathfrak p$.

{\bf\noindent Proof of Claim.} 
Since $\exp:\mathfrak {sl}_3(\mathbb R)\rightarrow SL_3(\mathbb R)$ is a local homeomorphism and $j$ is an injection, $(\exp\circ j)(\mathfrak p)$ is a submanifold of $G$ locally near the identity ${e}\in SL_3(\mathbb R)$.  Since $\mathfrak p\oplus \mathfrak a=\mathfrak {sl}_3(\mathbb R)$, where $\mathfrak a$ is the Lie algebra of $A$, we can conclude that $(\exp\circ j)(\mathfrak p)$ and $A$ intersects transversally at ${e}$. Since $\pi:SL_3(\mathbb R)\rightarrow SL_3(\mathbb R)/A$ is a fiber bundle which locally has a product structure,  $\mathcal P$ is a local homeomorphism between $\mathfrak p$ and $SL_3(\mathbb R)/A$.\,\,\,\,\,\,$\endpf$

Let $f$ be an arbitrary smooth function on $G/H$.  By Theorem \ref{HCiso}, we have 
\begin{equation*}
\begin{split}		&(D_{ij}f)(eA)=\left.\frac{\partial^2}{\partial x_{ij}\partial{x_{ji}}}f\left(\pi\left(\exp(x_{ij}X_{ij}+x_{ji}X_{ji})\right)\right)\right|_{x_{ij}=x_{ji}=0},\\		&(D_{ijk}f)(eA)=\left.\frac{\partial^3}{\partial x_{ij}\partial{x_{jk}}\partial x_{ki}}f\left(\pi\left(\exp(x_{ij}X_{ij}+x_{jk}X_{jk}+x_{ki}X_{ki})\right)\right)\right|_{x_{ij}=x_{jk}=x_{ki}=0},\\	\end{split}
\end{equation*}
for distinct $i,j,k\in\{1,2,3\}$. Then, in terms of the local coordinates $(x_{12},x_{13},x_{23},x_{21},x_{31},x_{32})$, we can conclude that the differential operators $D_{12}^kD_{23}^jD_{13}^iD_{123}^lD_{213}^m$, has leading terms at the origin $eA$ 
\begin{equation*}
\left(\frac{\partial}{\partial x_{12}}\right)^{k+l}\left(\frac{\partial}{\partial x_{23}}\right)^{j+l}\left(\frac{\partial}{\partial x_{31}}\right)^{l+i}\left(\frac{\partial}{\partial x_{21}}\right)^{k+m}\left(\frac{\partial}{\partial x_{13}}\right)^{m+i}\left(\frac{\partial}{\partial x_{32}}\right)^{j+m}.
\end{equation*} It is easy to verify that the differential operators appearing in (\ref{rel})
all have distinct leading terms at $eA$ . Then (\ref{rel}) holds if and only if all the coefficients $r_{kjilm}=0$. 

We complete the proof of Proposition \ref{basis}.\,\,\,$\endpf$
\smallskip

For convenience, we make the following convention. For a differential operator $D$, we write $D=O(M)$ if and only if the degree of $D$ is at most $M$. 

\begin{lemma}\label{KC}
The following equalities hold. 
\begin{equation*}
\begin{split}
&[D_{123},D_{12}^kD_{123}^lD_{213}^m]=kD_{12}^kD_{23}D_{123}^lD_{213}^m-kD_{12}^kD_{31}D_{123}^lD_{213}^m+O(2k+3l+3m+1),\\
&[D_{123},D_{23}^jD_{123}^lD_{213}^m]=-jD_{12}D_{23}^jD_{123}^lD_{213}^m+jD_{23}^jD_{31}D_{123}^lD_{213}^m+O(2j+3l+3m+1),\\
&[D_{123},D_{31}^iD_{123}^lD_{213}^m]=iD_{12}D_{31}^iD_{123}^lD_{213}^m-iD_{23}D_{31}^iD_{123}^lD_{213}^m+O(2i+3l+3m+1).\\
\end{split}
\end{equation*}
\end{lemma}
{\noindent\bf Proof of Lemma \ref{KC}.} We only prove the first identity in Lemma \ref{KC} here, as all the other ones can be proved by similar computation. By Proposition \ref{comm2}, we have 
\begin{equation*}
	\begin{split}
	&D_{123}D_{12}^kD_{123}^lD_{213}^m=D_{12}D_{123}D_{12}^{k-1}D_{123}^lD_{213}^m+(D_{23}D_{12}-D_{12}D_{31})D_{12}^{k-1}D_{123}^lD_{213}^m\\
	&=D_{12}^2D_{123}D_{12}^{k-2}D_{123}^lD_{213}^m+\sum\nolimits_{p=0}^1D_{12}^p(D_{23}D_{12}-D_{12}D_{31})D_{12}^{k-1-p}D_{123}^lD_{213}^m=\cdots\cdots\\
	&=D_{12}^kD_{123}^lD_{213}^mD_{123}+\sum\nolimits_{p=0}^{k-1}D_{12}^p(D_{23}D_{12}-D_{12}D_{31})D_{12}^{k-1-p}D_{123}^lD_{213}^m,\\
	\end{split}
\end{equation*}
where in the last step  (\ref{123213}) is applied. 

We complete the proof by counting the degree. $\endpf$
\medskip

Similarly, we have
\begin{lemma}\label{KJC}
The following equalities hold. 
\begin{equation*}
	\begin{split}
		&[D_{123},D_{12}^kD_{23}^jD_{123}^lD_{213}^m]=kD_{12}^kD_{23}^{j+1}D_{123}^lD_{213}^m+(j-k)D_{12}^{k-1}D_{23}^{j-1}D_{123}^{l+1}D_{213}^{m+1}\\			&\,\,\,\,\,\,\,\,\,\,\,\,\,\,\,\,\,\,\,\,\,\,-jD_{12}^{k+1}D_{23}^{j}D_{123}^lD_{213}^m+O(2k+2j+3l+3m+1),\\
		&[D_{123},D_{12}^kD_{31}^iD_{123}^lD_{213}^m]=(k-i)D_{12}^{k-1}D_{31}^{i-1}D_{123}^{l+1}D_{213}^{m+1}-kD_{12}^{k}D_{31}^{i+1}D_{123}^{l}D_{213}^{m}\\			&\,\,\,\,\,\,\,\,\,\,\,\,\,\,\,\,\,\,\,\,\,\,+iD_{12}^{k+1}D_{31}^{i}D_{123}^lD_{213}^m+O(2k+2i+3l+3m+1),\\
		&[D_{123},D_{23}^jD_{31}^iD_{123}^lD_{213}^m]=jD_{23}^{j}D_{31}^{i+1}D_{123}^lD_{213}^m+(i-j)D_{23}^{j-1}D_{31}^{i-1}D_{123}^{l+1}D_{213}^{m+1}\\			&\,\,\,\,\,\,\,\,\,\,\,\,\,\,\,\,\,\,\,\,\,\,-iD_{23}^{j+1}D_{31}^{i}D_{123}^{l}D_{213}^{m}+O(2j+2i+3l+3m+1).\\
	\end{split}
\end{equation*}
\end{lemma}

\begin{lemma}\label{KC1.1}
For $l,m\geq 0$, the following equalities hold. 
\begin{equation*}
	\begin{split}
		&[D_{12},D_{123}^lD_{213}^m]=-lD_{12}D_{23}D_{123}^{l-1}D_{213}^m+lD_{12}D_{31}D_{123}^{l-1}D_{213}^{m}\\
		&\,\,\,\,\,\,\,\,\,\,\,\,\,\,\,\,\,\,\,\,\,\,\,\,\,\,\,\,\,\,\,\,\,+mD_{12}D_{23}D_{123}^{l}D_{213}^{m-1}-mD_{12}D_{31}D_{123}^lD_{213}^{m-1}+O(3l+3m),\\
		&[D_{23},D_{123}^lD_{213}^m]=lD_{12}D_{23}D_{123}^{l-1}D_{213}^m-lD_{23}D_{31}D_{123}^{l-1}D_{213}^m\\
		&\,\,\,\,\,\,\,\,\,\,\,\,\,\,\,\,\,\,\,\,\,\,\,\,\,\,\,\,\,\,\,\,\,-mD_{12}D_{23}D_{123}^lD_{213}^{m-1}+mD_{23}D_{31}D_{123}^lD_{213}^{m-1}+O(3l+3m).\\
	\end{split}
\end{equation*}
\end{lemma}

\begin{lemma}\label{KC2}
For $k,j,i\geq 1$ and $l,m\geq 0$, the following equalities hold. 
\begin{equation*}
\begin{split}
&[D_{12},D_{23}^jD_{123}^lD_{213}^m]=-lD_{12}D_{23}^{j+1}D_{123}^{l-1}D_{213}^{m}+mD_{12}D_{23}^{j+1}D_{123}^{l}D_{213}^{m-1}\\
&\,\,\,\,\,\,\,\,\,\,\,\,\,\,\,\,\,\,\,\,\,\,\,\,\,\,\,\,\,\,\,\,\,+(j+l)D_{23}^{j-1}D_{123}^{l}D_{213}^{m+1}-(j+m)D_{23}^{j-1}D_{123}^{l+1}D_{213}^m+O(2j+3l+3m),\\
&[D_{12},D_{31}^iD_{123}^lD_{213}^m]=lD_{12}D_{31}^{i+1}D_{123}^{l-1}D_{213}^m-mD_{12}D_{31}^{i+1}D_{123}^lD_{213}^{m-1}\\
&\,\,\,\,\,\,\,\,\,\,\,\,\,\,\,\,\,\,\,\,\,\,\,\,\,\,\,\,\,\,\,\,\,+(i+m)D_{31}^{i-1}D_{123}^{l+1}D_{213}^m-(i+l)D_{31}^{i-1}D_{123}^{l}D_{213}^{m+1}+O(2i+3l+3m),\\
&[D_{23},D_{12}^kD_{123}^lD_{213}^m]=lD_{12}^{k+1}D_{23}D_{123}^{l-1}D_{213}^m-mD_{12}^{k+1}D_{23}D_{123}^lD_{213}^{m-1}\\
&\,\,\,\,\,\,\,\,\,\,\,\,\,\,\,\,\,\,\,\,\,\,\,\,\,\,\,\,\,\,\,\,\,-(k+l)D_{12}^{k-1}D_{123}^{l}D_{213}^{m+1}+(k+m)D_{12}^{k-1}D_{123}^{l+1}D_{213}^m+O(2k+3l+3m),\\
&[D_{23},D_{31}^iD_{123}^lD_{213}^m]=-lD_{23}D_{31}^{i+1}D_{123}^{l-1}D_{213}^m+mD_{23}D_{31}^{i+1}D_{123}^{l}D_{213}^{m-1}\\
&\,\,\,\,\,\,\,\,\,\,\,\,\,\,\,\,\,\,\,\,\,\,\,\,\,\,\,\,\,\,\,\,\,+(i+l)D_{31}^{i-1}D_{123}^{l}D_{213}^{m+1}-(i+m)D_{31}^{i-1}D_{123}^{l+1}D_{213}^{m}+O(2i+3l+3m).\\
\end{split}
\end{equation*}
\end{lemma}

\begin{lemma}\label{KJC2}
For $k,j,i\geq 1$ and $l,m\geq 0$, the following equalities hold. 
\begin{equation*}
\begin{split}
&[D_{12},D_{23}^jD_{31}^iD_{123}^lD_{213}^m]=-(i+l)D_{23}^{j}D_{31}^{i-1}D_{123}^{l}D_{213}^{m+1}+(i+m)D_{23}^{j} D_{31}^{i-1}D_{123}^{l+1}D_{213}^{m}\\			&\,\,\,\,\,\,\,\,+(j+l)D_{23}^{j-1} D_{31}^{i}D_{123}^lD_{213}^{m+1}-(j+m)D_{23}^{j-1}D_{31}^{i}D_{123}^{l+1}D_{213}^{m}+O(2j+2i+3l+3m),\\
&[D_{23},D_{12}^kD_{31}^iD_{123}^lD_{213}^m]=-(k+l)D_{12}^{k-1}D_{31}^{i}D_{123}^{l}D_{213}^{m+1}+(k+m)D_{12}^{k-1}D_{31}^{i}D_{123}^{l+1}D_{213}^{m}\\			&\,\,\,\,\,\,\,\,+(i+l)D_{12}^{k}D_{31}^{i-1} D_{123}^lD_{213}^{m+1}-(i+m)D_{12}^{k}D_{31}^{i-1} D_{123}^{l+1}D_{213}^{m}+O(2k+2i+3l+3m).\\
\end{split}
\end{equation*}  
\end{lemma}
{\noindent\bf Proof of Lemmas \ref{KJC}--\ref{KJC2}.} The proof is the same and we omit it for simplicity. \,\,\, $\endpf$

\begin{proposition}\label{center}
The center of $\mathbb D(SL_3(\mathbb R)/A)$ is a polynomial ring in $D_{123}+D_{213}$ and $D_{12}+D_{23}+D_{31}$.
\end{proposition}
{\bf\noindent Proof of Proposition \ref{center}.} 
By Propositions \ref{comm1}, \ref{comm2}, \ref{comm3}, it is easy to verify that $D_{123}+D_{213}$, $D_{12}+D_{23}+D_{31}$ are in the center.
We will use the idea of \cite{HY} to solve a large system of linear equations inductively.

We claim that the algebra generated by $D_{123}+D_{213}$ and $D_{12}+D_{23}+D_{31}$ is a polynomial ring. Otherwise, there is a nontrivial relation 
\begin{equation*}
0=\sum\nolimits_{\substack{K,M\in\mathbb Z{^{\geq 0}}}}C_{KM}\cdot(D_{12}+D_{23}+D_{31})^K(D_{123}+D_{213})^M=:R.
\end{equation*}
Define $A:=\max\{\,K\,|\,C_{KM}\neq 0\,\}$ and $B:=\max\{\,M\,|\,C_{AM}\neq 0\,\}$. Expanding $R$ in terms of the basis in (\ref{47}), we have  
\begin{equation*}
	\begin{split}
		R=&\sum_{\substack{l,m\in\mathbb Z{^{\geq 0}}}}r_{000lm}D_{123}^lD_{213}^m+\sum_{\substack{k\in\mathbb Z^+\\l,m\in\mathbb Z{^{\geq 0}}}}r_{k00lm}D_{12}^kD_{123}^lD_{213}^m+\sum_{\substack{j\in\mathbb Z^+\\l,m\in\mathbb Z{^{\geq 0}}}}r_{0j0lm}D_{23}^jD_{123}^lD_{213}^m\\
		&+\sum_{\substack{i\in\mathbb Z^+\\l,m\in\mathbb Z{^{\geq 0}}}}r_{00ilm}D_{31}^iD_{123}^lD_{213}^m+	
		\sum_{\substack{k,j\in\mathbb Z^+\\l,m\in\mathbb Z{^{\geq 0}}}}r_{kj0lm} D_{12}^kD_{23}^jD_{123}^lD_{213}^m\\
		&+\sum_{\substack{k,i\in\mathbb Z^+\\l,m\in\mathbb Z{^{\geq 0}}}}r_{k0ilm} D_{12}^kD_{31}^iD_{123}^lD_{213}^m+\sum_{\substack{j,i\in\mathbb Z^+\\l,m\in\mathbb Z{^{\geq 0}}}}r_{0jilm}D_{23}^jD_{31}^iD_{123}^lD_{213}^m.\\
	\end{split}
\end{equation*}
It is easy to verify that $r_{A000B}=C_{AB}\neq 0$, which is a contradiction.

In what follows, we will prove that the center of $\mathbb D(SL_3(\mathbb R)/A)$ is generated by $D_{123}+D_{213}$ and $D_{12}+D_{23}+D_{31}$. By Proposition \ref{basis}, write any $z\in \mathbb D(SL_3(\mathbb R)/A)$ as
\begin{equation*}
	\begin{split}
		z=&\sum_{\substack{l,m\in\mathbb Z{^{\geq 0}}}}a_{000lm}D_{123}^lD_{213}^m+\sum_{\substack{k\in\mathbb Z^+\\l,m\in\mathbb Z{^{\geq 0}}}}a_{k00lm}D_{12}^kD_{123}^lD_{213}^m+\sum_{\substack{j\in\mathbb Z^+\\l,m\in\mathbb Z{^{\geq 0}}}}a_{0j0lm}D_{23}^jD_{123}^lD_{213}^m\\
		&+\sum_{\substack{i\in\mathbb Z^+\\l,m\in\mathbb Z{^{\geq 0}}}}a_{00ilm}D_{31}^iD_{123}^lD_{213}^m+	
		\sum_{\substack{k,j\in\mathbb Z^+\\l,m\in\mathbb Z{^{\geq 0}}}}a_{kj0lm} D_{12}^kD_{23}^jD_{123}^lD_{213}^m\\
		&+\sum_{\substack{k,i\in\mathbb Z^+\\l,m\in\mathbb Z{^{\geq 0}}}}a_{k0ilm} D_{12}^kD_{31}^iD_{123}^lD_{213}^m+\sum_{\substack{j,i\in\mathbb Z^+\\l,m\in\mathbb Z{^{\geq 0}}}}a_{0jilm}D_{23}^jD_{31}^iD_{123}^lD_{213}^m.\\
	\end{split}
\end{equation*}
We make the convention that $a_{KJILM}=0$ when $K,J,I,L$, or $M$ is strictly negative.
Applying Lemmas \ref{KC}, \ref{KJC}, we have 
\small
\begin{equation*}
\begin{split}
0=&[D_{123},z]=\sum_{\substack{k\in\mathbb Z^+,\,\,l,m\in\mathbb Z{^{\geq 0}}}}a_{k00lm}\left(kD_{12}^kD_{23}D_{123}^lD_{213}^m-kD_{12}^kD_{31}D_{123}^lD_{213}^m\right)\,\,\,\,\,\,\,\,\,\,\,\,\,\,\,\,\,\,\,\,\,\,\,\,\,\,\,\,\,\,\,\,\,\,\,\,\,\,\,\,\,\,\,\,\,\,\,\,\,\,\,\,\,\,\,\,\,\,\,\,\,\,\,\,\,\,\\
&+\sum_{\substack{j\in\mathbb Z^+,\,\,l,m\in\mathbb Z{^{\geq 0}}}}a_{0j0lm}\left(-jD_{12}D_{23}^jD_{123}^lD_{213}^m+jD_{23}^jD_{31}D_{123}^lD_{213}^m\right)\\
&+\sum_{\substack{i\in\mathbb Z^+,\,\,l,m\in\mathbb Z{^{\geq 0}}}}a_{00ilm}\left(iD_{12}D_{31}^iD_{123}^lD_{213}^m-iD_{23}D_{31}^iD_{123}^lD_{213}^m\right)\\
\end{split}
\end{equation*}
\normalsize
\vspace{-.12in}
\footnotesize
\begin{equation*}
	\begin{split}
	\,\,	&+\sum_{\substack{k,j\in\mathbb Z^+,\,\,l,m\in\mathbb Z{^{\geq 0}}}}a_{kj0lm}\left(kD_{12}^kD_{23}^{j+1}D_{123}^lD_{213}^m+(j-k)D_{12}^{k-1}D_{23}^{j-1}D_{123}^{l+1}D_{213}^{m+1}-jD_{12}^{k+1}D_{23}^{j}D_{123}^lD_{213}^m\right)\\
		&+\sum_{\substack{k,i\in\mathbb Z^+,\,\,l,m\in\mathbb Z{^{\geq 0}}}}a_{k0ilm}\left((k-i)D_{12}^{k-1}D_{31}^{i-1}D_{123}^{l+1}D_{213}^{m+1}-kD_{12}^{k}D_{31}^{i+1}D_{123}^{l}D_{213}^{m}+iD_{12}^{k+1}D_{31}^{i}D_{123}^lD_{213}^m\right)\\
		&+\sum_{\substack{j,i\in\mathbb Z^+,\,\,l,m\in\mathbb Z{^{\geq 0}}}}a_{0jilm}\left(jD_{23}^{j}D_{31}^{i+1}D_{123}^lD_{213}^m+(i-j)D_{23}^{j-1}D_{31}^{i-1}D_{123}^{l+1}D_{213}^{m+1}-iD_{23}^{j+1}D_{31}^{i}D_{123}^{l}D_{213}^{m}\right).\\		
	\end{split}
\end{equation*}
\normalsize
By the linear independence, we can conclude that, for $L,M\geq 0$ and $K,J\geq 1$, 
\begin{equation}\label{123KJ}
	Ka_{K(J-1)0LM}+(J-K)a_{(K+1)(J+1)0(L-1)(M-1)}-Ja_{(K-1)J0LM}=0;
\end{equation}
for $L,M\geq 0$ and $K,I\geq 1$,
\begin{equation}\label{79}
	(K-I)a_{(K+1)0(I+1)(L-1)(M-1)}-Ka_{K0(I-1)LM}+Ia_{(K-1)0ILM}=0;
\end{equation}
for $L,M\geq 0$ and $J,I\geq 1$,
\begin{equation}\label{80}Ja_{0J(I-1)LM}+(I-J)a_{0(J+1)(I+1)(L-1)(M-1)}-Ia_{0(J-1)ILM}=0.
\end{equation}

We have by Lemmas \ref{KC2}, \ref{KJC2} that
\small
\begin{equation*}
	\begin{split}
		0&=[D_{12},z]=\sum_{\substack{l,m\in\mathbb Z{^{\geq 0}}}}a_{000lm}\left(-lD_{12}D_{23}D_{123}^{l-1}D_{213}^m+lD_{12}D_{31}D_{123}^{l-1}D_{213}^{m}+mD_{12}D_{23}D_{123}^{l}D_{213}^{m-1}\right.\\
		&-m\left.D_{12}D_{31}D_{123}^lD_{213}^{m-1}\right)+\sum_{\substack{k\in\mathbb Z^+,\,\,l,m\in\mathbb Z{^{\geq 0}}}}a_{k00lm}\left(-lD_{12}^{k+1}D_{23}D_{123}^{l-1}D_{213}^m+lD_{12}^{k+1}D_{31}D_{123}^{l-1}D_{213}^{m}\right.\\
		&+mD_{12}^{k+1}D_{23}D_{123}^{l}D_{213}^{m-1}\left.-mD_{12}^{k+1}D_{31}D_{123}^lD_{213}^{m-1}\right)+\sum_{\substack{j\in\mathbb Z^+,\,\,l,m\in\mathbb Z{^{\geq 0}}}}a_{0j0lm}\left(-lD_{12}D_{23}^{j+1}D_{123}^{l-1}D_{213}^{m}\right.\,\,\,\\
		&+mD_{12}D_{23}^{j+1}D_{123}^{l}D_{213}^{m-1}\left.+(j+l)D_{23}^{j-1}D_{123}^{l}D_{213}^{m+1}-(j+m)D_{23}^{j-1}D_{123}^{l+1}D_{213}^m\right)\\
		&+\sum_{\substack{i\in\mathbb Z^+,\,\,l,m\in\mathbb Z{^{\geq 0}}}}a_{00ilm}\left(lD_{12}D_{31}^{i+1}D_{123}^{l-1}D_{213}^m-mD_{12}D_{31}^{i+1}D_{123}^lD_{213}^{m-1}-(i+l)D_{31}^{i-1}D_{123}^{l}D_{213}^{m+1}\right.\\
	\end{split}
\end{equation*}
\normalsize
\vspace{-.13in}
\small
\begin{equation*}
	\begin{split}
		\,\,	&\left.+(i+m)D_{31}^{i-1}D_{123}^{l+1}D_{213}^m\right)+\sum_{\substack{k,j\in\mathbb Z^+,\,\,l,m\in\mathbb Z{^{\geq 0}}}}a_{kj0lm}\left(-lD_{12}^{k+1} D_{23}^{j+1} D_{123}^{l-1}D_{213}^{m}\right.\\
		&\left.+mD_{12}^{k+1}D_{23}^{j+1}D_{123}^{l} D_{213}^{m-1}+(j+l)D_{12}^k D_{23}^{j-1} D_{123}^lD_{213}^{m+1}-(j+m)D_{12}^{k}D_{23}^{j-1}D_{123}^{l+1} D_{213}^{m}\right)\\
		&+\sum_{k,i\in\mathbb Z^+,\,\,l,m\in\mathbb Z{^{\geq 0}}}a_{k0ilm}
		\left(lD_{12}^{k+1}D_{31}^{i+1}D_{123}^{l-1}D_{213}^m-mD_{12}^{k+1}D_{31}^{i+1}D_{123}^lD_{213}^{m-1}-(i+l)D_{12}^{k}D_{31}^{i-1} D_{123}^{l}D_{213}^{m+1}\right.\\
		&\left.+(i+m)D_{12}^{k} D_{31}^{i-1}D_{123}^{l+1}D_{213}^{m}\right)+\sum_{\substack{j,i\in\mathbb Z^+,\,\,l,m\in\mathbb Z{^{\geq 0}}}}a_{0jilm}\left(-(i+l)D_{23}^{j}D_{31}^{i-1}D_{123}^{l}D_{213}^{m+1}\right.\\
		&\left.+(i+m)D_{23}^{j} D_{31}^{i-1}D_{123}^{l+1}D_{213}^{m}+(j+l)D_{23}^{j-1} D_{31}^{i}D_{123}^lD_{213}^{m+1}-(j+m)D_{23}^{j-1}D_{31}^{i}D_{123}^{l+1}D_{213}^{m}\right).\\		
	\end{split}
\end{equation*}
\normalsize
Then 
for $L,M\geq 0$ and $K,I\geq 1$,
\begin{equation}\label{86}
\begin{split}
&(L+1)a_{(K-1)0(I-1)(L+1)M}-(M+1)a_{(K-1)0(I-1)L(M+1)}\,\,\,\,\,\,\\
&-(I+L+1)a_{K0(I+1)L(M-1)}+(I+M+1)a_{K0(I+1)(L-1)M}=0.
\end{split}
\end{equation}

Similarly,
\small\begin{equation*}
	\begin{split}
		&0=[D_{23},z]=\sum_{\substack{l,m\in\mathbb Z{^{\geq 0}}}}a_{000lm}\left(lD_{12}D_{23}D_{123}^{l-1}D_{213}^m-lD_{23}D_{31}D_{123}^{l-1}D_{213}^m-mD_{12}D_{23}D_{123}^lD_{213}^{m-1}\right.\\
		&+m\left.D_{23}D_{31}D_{123}^lD_{213}^{m-1}\right)+\sum_{\substack{k\in\mathbb Z^+,\,\,l,m\in\mathbb Z{^{\geq 0}}}}a_{k00lm}\left(lD_{12}^{k+1}D_{23}D_{123}^{l-1}D_{213}^m-mD_{12}^{k+1}D_{23}D_{123}^lD_{213}^{m-1}\right.\\
		&-(k+l)D_{12}^{k-1}D_{123}^{l}D_{213}^{m+1}\left.+(k+m)D_{12}^{k-1}D_{123}^{l+1}D_{213}^m\right)+\sum_{\substack{j\in\mathbb Z^+,\,\,l,m\in\mathbb Z{^{\geq 0}}}}a_{0j0lm}\left(lD_{12}D_{23}^{j+1}D_{123}^{l-1}D_{213}^m\right.\\
		&-lD_{23}^{j+1}D_{31}D_{123}^{l-1}D_{213}^m\left.-mD_{12}D_{23}^{j+1}D_{123}^lD_{213}^{m-1}+mD_{23}^{j+1}D_{31}D_{123}^lD_{213}^{m-1}\right)\\
		&+\sum_{\substack{i\in\mathbb Z^+,\,\,l,m\in\mathbb Z{^{\geq 0}}}}a_{00ilm}\left(-lD_{23}D_{31}^{i+1}D_{123}^{l-1}D_{213}^m+mD_{23}D_{31}^{i+1}D_{123}^{l}D_{213}^{m-1}+(i+l)D_{31}^{i-1}D_{123}^{l}D_{213}^{m+1}\right.\\
		&-(i+m)\left.D_{31}^{i-1}D_{123}^{l+1}D_{213}^{m}\right)+\sum_{\substack{k,j\in\mathbb Z^+,\,\,l,m\in\mathbb Z{^{\geq 0}}}}a_{kj0lm}\left(lD_{12}^{k+1}D_{23}^{j+1}D_{123}^{l-1}D_{213}^m-mD_{12}^{k+1}D_{23}^{j+1}D_{123}^lD_{213}^{m-1}\right.\,\,\,\,\,\,\,\\
	\end{split}
\end{equation*}
\normalsize
\vspace{-.1in}
\small
\begin{equation*}
	\begin{split}
		&-(k+l)D_{12}^{k-1}D_{23}^{j}D_{123}^{l}D_{213}^{m+1}\left.+(k+m)D_{12}^{k-1}D_{23}^{j}D_{123}^{l+1}D_{213}^m\right)\\
		&+\sum_{\substack{k,i\in\mathbb Z^+,\,\,l,m\in\mathbb Z{^{\geq 0}}}}a_{k0ilm}\left(-(k+l)D_{12}^{k-1}D_{31}^{i}D_{123}^{l}D_{213}^{m+1}+(k+m)D_{12}^{k-1}D_{31}^{i}D_{123}^{l+1}D_{213}^{m}\right.\\
		&+(i+l)D_{12}^{k}D_{31}^{i-1} D_{123}^lD_{213}^{m+1}\left.-(i+m)D_{12}^{k}D_{31}^{i-1} D_{123}^{l+1}D_{213}^{m}\right)+\sum_{\substack{j,i\in\mathbb Z^+,\,\,l,m\in\mathbb Z{^{\geq 0}}}}a_{0jilm}\left(-lD_{23}^{j+1}D_{31}^{i+1}D_{123}^{l-1} D_{213}^{m}\right.\,\,\\
		&+mD_{23}^{j+1}D_{31}^{i+1}D_{123}^{l}D_{213}^{m-1}\left.+(i+l)D_{23}^{j}D_{31}^{i-1}D_{123}^lD_{213}^{m+1}-(i+m)D_{23}^{j}D_{31}^{i-1} D_{123}^{l+1}D_{213}^{m}\right).\\		
	\end{split}
\end{equation*}
\normalsize
Then for $L,M\geq 0$ and $K,J\geq 1$, 
\begin{equation}\label{23KJ}
	\begin{split}
		&(L+1)a_{(K-1)(J-1)0(L+1)M}-(M+1)a_{(K-1)(J-1)0L(M+1)}\\
		&-(K+L+1)a_{(K+1)J0L(M-1)}+(K+M+1)a_{(K+1)J0(L-1)M}=0,
	\end{split}
\end{equation}
and for $L,M\geq 0$ and $J,I\geq 1$,
\begin{equation}\label{92}
\begin{split}
&-(L+1)a_{0(J-1)(I-1)(L+1)M}+(M+1)a_{0(J-1)(I-1)L(M+1)}\,\,\,\,\,\,\\
&+(I+L+1)a_{0J(I+1)L(M-1)}-(I+M+1)a_{0J(I+1)(L-1)M}=0.
\end{split}
\end{equation}

\noindent{\bf Claim I.} If $a_{K000M}=0$ for $K,M\geq 0$, then $a_{KJ0LM}\equiv 0$ for $K,J,L,M\geq 0$.
\smallskip

\noindent{\bf Proof of Claim I.} We prove  by induction on $L$. Setting $L=0$ in (\ref{123KJ}), we have for $M\geq0$ and $K,J\geq 1$ that $Ka_{K(J-1)00M}=Ja_{(K-1)J00M}$. Hence, when $M\geq0$ and $K,J\geq 1$,
\begin{equation*}
a_{(K-1)J00M}=\frac{(K+J-1)!}{J!\cdot(K-1)!}\cdot a_{(K+J-1)000M}=0.
\end{equation*}
We thus prove Claim I when $L=0$.

Suppose Claim I holds for $0\leq L\leq l$. By setting $L=l$ in (\ref{23KJ}), we can conclude that for $M\geq 0$ and $K,J\geq 1$, $a_{(K-1)(J-1)0(l+1)M}=0$.
Then, Claim I holds for $L=l+1$. Therefore, we complete the proof of Claim I.   \,\,\,\,$\endpf$
\smallskip

\noindent{\bf Claim II.} If $a_{K000M}=0$ for $K,M\geq 0$, then $a_{K0ILM}\equiv 0$ for $K,I,L,M\geq 0$.
\smallskip

\noindent{\bf Proof of Claim II.} Similarly, we prove by induction on $L$. Setting $L=0$ in (\ref{79}), we gave
$Ka_{K0(I-1)0M}=Ia_{(K-1)0I0M}$,  for $M\geq0,K\geq 1,I\geq 1$,
and hence when $M\geq0, K\geq 1,I\geq 1,$
\begin{equation*}
	a_{(K-1)0I0M}=\frac{(K+I-1)!}{I!\cdot(K-1)!}\cdot a_{(K+I-1)000M}=0,
\end{equation*}
Suppose Claim II holds for $0\leq L\leq l$. Setting $L=l$ in (\ref{86}), we can conclude that
$a_{(K-1)0(I-1)(l+1)M}=0$, for $M\geq 0,K\geq 1,I\geq 1$.

We complete the proof of Claim II.   \,\,\,\,$\endpf$
\smallskip

\noindent{\bf Claim III.} If $a_{K000M}=0$ for $K,M\geq 0$, then $a_{0JILM}\equiv 0$ for $K,J,L,M\geq 0$.
\smallskip

\noindent{\bf Proof of Claim III.} By Claim II, it is clear that for $I,M\geq 0$, $a_{00I0M}=0$. Similarly to the proofs of Claims I and II, we can then complete the proof of Claim III by (\ref{80}) and (\ref{92}). \,\,\,\,$\endpf$
\smallskip

Define
\begin{equation*}
	z^{\prime}:=z-\sum\nolimits_{\substack{K,M\in\mathbb Z{^{\geq 0}}}}a_{K000M}(D_{12}+D_{23}+D_{31})^K(D_{123}+D_{213})^M.
\end{equation*}
Then $z^{\prime}$ belongs to the center and can be written as
\begin{equation*}
	\begin{split}
		z^{\prime}=&\sum_{\substack{l,m\in\mathbb Z{^{\geq 0}}}}a_{000lm}^{\prime}D_{123}^lD_{213}^m+\sum_{\substack{k\in\mathbb Z^+\\l,m\in\mathbb Z{^{\geq 0}}}}a_{k00lm}^{\prime}D_{12}^kD_{123}^lD_{213}^m+\sum_{\substack{j\in\mathbb Z^+\\l,m\in\mathbb Z{^{\geq 0}}}}a_{0j0lm}^{\prime}D_{23}^jD_{123}^lD_{213}^m\\
		&+\sum_{\substack{i\in\mathbb Z^+\\l,m\in\mathbb Z{^{\geq 0}}}}a_{00ilm}^{\prime}D_{31}^iD_{123}^lD_{213}^m+	
		\sum_{\substack{k,j\in\mathbb Z^+\\l,m\in\mathbb Z{^{\geq 0}}}}a_{kj0lm}^{\prime} D_{12}^kD_{23}^jD_{123}^lD_{213}^m\\
		&+\sum_{\substack{k,i\in\mathbb Z^+\\l,m\in\mathbb Z{^{\geq 0}}}}a_{k0ilm}^{\prime} D_{12}^kD_{31}^iD_{123}^lD_{213}^m+\sum_{\substack{j,i\in\mathbb Z^+\\l,m\in\mathbb Z{^{\geq 0}}}}a_{0jilm}^{\prime}D_{23}^jD_{31}^iD_{123}^lD_{213}^m.\\
	\end{split}
\end{equation*}
such that for $K,M\geq 0$, $a^{\prime}_{K000M}=0$. Now by Claims I, II, III, we can conclude that $z^{\prime}=0$.

We complete the proof of Theorem \ref{center}. \,\,\,\,$\endpf$
\medskip

{\bf\noindent Proof of Theorem \ref{main1}.}
It follows from Lemma \ref{generators}, and Propositions \ref{comm1}, \ref{comm2}, \ref{comm4}, \ref{comm3},  \ref{basis}, \ref{center}.
\,\,\,\,$\endpf$


\section{Essential Self-Adjointness}\label{sellf}

\subsection{Essential self-adjointness of the central elements}\label{central}

A classical result due to \cite{Seg,NS} says that if $\tau$ is the quasi-regular representation of $G$ in $L^2(G/H)$,  $D$ is a symmetric element in the center of the universal enveloping algebra $U(\mathfrak{g})$ of the Lie algebra $\mathfrak{g}$ of $G$, and the operator $\tau(D)$ is an invariant differential operator on $X$, then $\tau(D)$ is essentially self-adjoint. Letting $\tau$ be the natural representation of $G$ in $L^2(G/H)$ (i.e. by translation), for any element $D$ in the center of $U(\mathfrak{g})$, $\tau(D)$ equals $\mu(R(D))$ as defined in \eqref{mu}. 
\medskip

{\bf\noindent Proof of Proposition \ref{easy}.} Recall that the center of $U(\mathfrak{sl}_3(\mathbb R))$ is a polynomial 
ring $Z=\mathbb R[h, k]$ (see Page 984 of \cite{Ca}), where
\begin{equation*}
\begin{split}
&h=-(X_{11}-X_{22})^2+(X_{11}-X_{22})X_{11}-X_{11}^2-3E_{12}E_{21}-3E_{13}E_{31}-3E_{23}E_{32}+3X_{11}, \\
&k=2(X_{11}-X_{22})^3-3(X_{11}-X_{22})^2X_{11}-3(X_{11}-X_{22})X_{11}^2+2X_{11}^3\\
&\,\,\,\,\,\,\,\,\,+9E_{12}E_{21}(X_{11}-X_{22})-18E_{12}E_{21}X_{11}-18E_{13}E_{31}(X_{11}-X_{22})+9E_{13}E_{31}X_{11}\,\,\,\,\,\,\,\,\,\,\,\,\\
&\,\,\,\,\,\,\,\,\,+9E_{23}E_{32}(X_{11}-X_{22})+9E_{23}E_{32}X_{11}-27E_{12}E_{23}E_{31}-27E_{21}E_{13}E_{32}\\
&\,\,\,\,\,\,\,\,\,+18(X_{11}-X_{22})X_{11}-9X_{11}^2-18(X_{11}-X_{22})+9X_{11}.    
\end{split}    
\end{equation*}
It is clear that $R(h),R(k)\in\mathbb D^A(SL_3(\mathbb R))$, and
\begin{equation*}
\begin{split}
&\mu(R(h))=-3\mu(R(E_{12}E_{21}))-3\mu(R(E_{13}E_{31}))-3\mu(R(E_{23}E_{32}))=-3(D_{12}+D_{13}+D_{23}),\\
&\mu(R(k))=-27\mu(R(E_{12}E_{31}E_{23}))-27\mu(R(E_{21}E_{13}E_{32}))=-27(D_{123}+D_{213}).\\
\end{split}    
\end{equation*}
By the above result due to \cite{Seg,NS}, we complete the proof of Proposition \ref{easy}.

\subsection{Coordinate charts induced by the Euler angles}\label{coordinate}
In this subsection, we describe certain  geometric properties of $X=SL_3(\mathbb R)/A$ for explicit computation. For convenience, in the following, we will use $X$ and $SL_3(\mathbb R)/A$ interchangeably.

According to the Iwasawa decomposition $SL_3(\mathbb R)=KNA$,  each $x\in SL_3(\mathbb R)$ can be uniquely written as $x=k_xn_xa_x$, where $k_x$ is an unimodular orthogonal matrix, $n_x$ is an upper unitriangular matrix, and $a_x$ is a unimodular diagonal matrix with positive entries; more explicitly, 
\begin{equation}\label{iwasawa}
x=\left(\begin{matrix}k_{11}&k_{12}&k_{13}\\		k_{21}&k_{22}&k_{23}\\		k_{31}&k_{32}&k_{33}\end{matrix}\right)\left(\begin{matrix}1&z_{12}&z_{13}\\		0&1&z_{23}\\
0&0&1\\\end{matrix}\right)\left(\begin{matrix}\lambda_1&0&0\\	0&\lambda_{2}&0\\
0&0&\lambda_{1}^{-1}\lambda_2^{-1}\end{matrix}\right).
\end{equation}
Since $SL_3(\mathbb R)/A$ is diffeomorphic to $SO_3(\mathbb R) \times\mathbb R^3$, without loss of generality, we can identify each $x\in SL_3(\mathbb R)/A$ with $(k_x,z_{12},z_{13},z_{23})\in SO_3(\mathbb R) \times\mathbb R^3$.

We parametrize $SO_3(\mathbb R)$ by the following variant of Euler angles (Tait-Bryan angles). Define $\Gamma:[-\pi,\pi]\times[-\frac{\pi}{2},\frac{\pi}{2}]\times[-\pi,\pi]\rightarrow SO_3(\mathbb R)$ by $(\alpha,\beta,\gamma)\mapsto k_x$, where
\footnotesize
\begin{equation}\label{tb}
k_x:=\left(\begin{matrix}k_{11}&k_{12}&k_{13}\\		k_{21}&k_{22}&k_{23}\\		k_{31}&k_{32}&k_{33}\end{matrix}\right):=\left(\begin{matrix}\cos{\beta}\cos{\gamma}&-\sin{\beta}&\cos{\beta}\sin{\gamma}\\\sin{\alpha}\sin{\gamma}+\cos{\alpha}\cos{\gamma}\sin{\beta}&\cos{\alpha}\cos{\beta}&\cos{\alpha}\sin{\beta}\sin{\gamma}-\cos{\gamma}\sin{\alpha}\\\cos{\gamma}\sin{\alpha}\sin{\beta}-\cos{\alpha}\sin{\gamma}&\cos{\beta}\sin{\alpha}&\cos{\alpha}\cos{\gamma}+\sin{\alpha}\sin{\beta}\sin{\gamma}\end{matrix}\right).
\end{equation}
\normalsize
It is easy to verify that when $-\frac{\pi}{2}<\alpha<\frac{\pi}{2}$, $-\frac{\pi}{3}<\beta< \frac{\pi}{3}$,  $-\frac{\pi}{2}<\gamma<\frac{\pi}{2}$, $\Gamma$ is bijective onto the image, and the inverse map is given by
\begin{equation}\label{inverse}
\alpha =\arctan \left({\frac {k_{32}}{k_{22}}}\right),\,\,\beta=\arctan \left(\frac{-k_{12}}{\sqrt{1-k_{12}^2}}\right),\,\,\gamma=\arctan \left({\frac {k_{13}}{k_{11}}}\right).
\end{equation}
\begin{definition}\label{u0}  We denote by $\underline{U}_{EI}$ the image under $\Gamma$ of  $-\frac{\pi}{2}<\alpha<\frac{\pi}{2}$, $-\frac{\pi}{3}<\beta< \frac{\pi}{3}$,  $-\frac{\pi}{2}<\gamma<\frac{\pi}{2}$. It is clear that $\underline{U}_{EI}$ is a neighborhood of the identity element in $SO_3(\mathbb R)$. 
\end{definition}

Denote by $(k_0\cdot\underline{U}_{EI})$ the left translation of $\underline{U}_{EI}$ by $k_0$. Then,  we can construct an atlas $\left\{\left(\left(k_0\cdot\underline{U}_{EI}\right)\times\mathbb R^3,\Lambda_{k_0}\right)\right\}_{k_0\in SO_3(\mathbb R)}$ of $SL_3(\mathbb R)/A$ as follows. For each $k_0\in SO_3(\mathbb R)$, define a map $\Lambda_{k_0}:\left(k_0\cdot\underline{U}_{EI}\right)\times\mathbb R^3\rightarrow(-\frac{\pi}{2},\frac{\pi}{2})\times(-\frac{\pi}{3},\frac{\pi}{3})\times(-\frac{\pi}{2},\frac{\pi}{2})\times\mathbb R^3$ by
\begin{equation*}
\left(k_0\cdot \left(\begin{matrix}k_{11}&k_{12}&k_{13}\\		k_{21}&k_{22}&k_{23}\\		k_{31}&k_{32}&k_{33}\end{matrix}\right),z_{12},z_{13},z_{23}\right)\mapsto (\alpha,\beta,\gamma,z_{12},z_{13},z_{23}),\\
\end{equation*}
where $\alpha,\beta,\gamma$ are related to $(k_{ij})$ by (\ref{inverse}).
\begin{definition}\label{ei}
We call the coordinates $(\alpha,\beta,\gamma,z_{12},z_{13},z_{23})$, which is associated with the coordinate chart $\left(\left(k_0\cdot\underline{U}_{EI}\right)\times\mathbb R^3,\Lambda_{k_0}\right)$ as above, the Euler-Iwasawa coordinates of $SL_3(\mathbb R)/A$ attached to $k_0$.
\end{definition}

Notice that the Euler-Iwasawa coordinates attached to $k_0$, combined with the global coordinates $(\lambda_1,\lambda_2)$ of $A$ (see (\ref{iwasawa})), parametrize the open neighborhood $\left(k_0\cdot\underline{U}_{EI}\right)\times\mathbb R^3\times \mathbb R^2\subset SL_3(\mathbb R)$. We denote the corresponding atlas of $SL_3(\mathbb R)$ by $\left\{\left(\left(k_0\cdot\underline{U}_{EI}\right)\times\mathbb R^3\times \mathbb R^2,\widehat\Lambda_{k_0}\right)\right\}_{k_0\in SO_3(\mathbb R)}$. 
\begin{definition}\label{eisl}
We call the coordinates $(\alpha,\beta,\gamma,z_{12}$, $z_{13},z_{23},\lambda_1,\lambda_2)$, which is associated with the coordinate chart $\left(\left(k_0\cdot\underline{U}_{EI}\right)\times\mathbb R^3\times \mathbb R^2,\widehat\Lambda_{k_0}\right)$ as above,  the Euler-Iwasawa coordinates of $SL_3(\mathbb R)$ attached to $k_0$.
\end{definition}

One can verify that an invariant measure on $X$ in the coordinate chart $\left(\left(k_0\cdot\underline{U}_{EI}\right)\times\mathbb R^3,\Lambda_{k_0}\right)$ takes the form 
\begin{equation}\label{hm}
{\rm d}\mu={\rm d}k\  {\rm d}n=\cos\beta\, {\rm d}\alpha\  {\rm d}\beta\  {\rm d}\gamma\  {\rm d}z_{12}\  {\rm d}z_{13}\  {\rm d}	z_{23}.	
\end{equation}

\subsection{Represent \texorpdfstring{$D_{12}$}{dd} as left derivatives}\label{prop}

We first recall the notion of left derivatives as in \cite{Ba}. 
Notice that the infinitesimal action $L$ on $C^{\infty}(SL_3(\mathbb R))$ induced by the left regular representation, maps $U(\mathfrak {sl}_3(\mathbb R))$ into the algebra of differentials on $SL_3(\mathbb R)$. More precisely, for $u=\sum X_{1}X_2\cdots X_k\in U(\mathfrak {sl}_3(\mathbb R))$,  and $f \in C^{\infty}(SL_3(\mathbb R))$,
\begin{equation*}
\left(L(u) f\right)(g):=(L_u f)(g):=\sum \left.\frac{\partial}{\partial t_1}\right|_{t_1=0}\cdots\left.\frac{\partial}{\partial t_k}\right|_{t_k=0}f(\exp \left(-t_kX_k\right)\cdots\exp \left(-t_1X_1\right)g).
\end{equation*}
Then, $L$ also induces differential operators on $SL_3(\mathbb R)/A$, that is, for $u=\sum X_{1}X_2\cdots X_k\in U(\mathfrak {sl}_3(\mathbb R))$, and $f \in C^{\infty}(SL_3(\mathbb R)/A)$,
\begin{equation*}
\left(L(u) f\right)(xA):=\sum \left.\frac{\partial}{\partial t_1}\right|_{t_1=0}\cdots\left.\frac{\partial}{\partial t_k}\right|_{t_k=0}f(\exp \left(-t_kX_k\right)\cdots\exp \left(-t_1X_1\right)xA).
\end{equation*}

\begin{definition}\label{ld}
We call a differential operator on $SL_3(\mathbb R)/A$ a left derivative if it is induced by an element $u\in U(\mathfrak{sl}_3(\mathbb R))$ through the above infinitesimal action $L$. For convenience, denote left derivatives by $L_u$, $u\in U(\mathfrak{sl}_3(\mathbb R))$, as well. 
\end{definition}


\begin{lemma}\label{eid12} 
In each coordinate chart $\left(\left(k_0\cdot\underline{U}_{EI}\right)\times\mathbb R^3,\Lambda_{k_0}\right)$,  $D_{12}$ takes the  form 
\begin{equation}\label{feid12}
	\begin{split}
		&D_{12}=\left(\sec\beta\sin\gamma\right)\frac{\partial^2 }{\partial \alpha\partial z_{12}}+\left(\cos\gamma\right)\frac{\partial^2}{\partial \beta\partial z_{12}} +\left(\tan\beta\sin\gamma\right) \frac{\partial^2}{\partial \gamma\partial z_{12}}\\
		&\,\,\,\,\,\,\,\,\,\,\,\,\,+\left(z_{12}^2+1\right)\frac{\partial^2}{\partial z_{12}\partial z_{12}}+\left(z_{23}\right)\frac{\partial^2}{\partial z_{13}\partial z_{12}}+\left(-z_{13}\right)\frac{\partial^2}{\partial z_{23}\partial z_{12}}+2z_{12}\frac{\partial }{\partial z_{12}},\\
	\end{split}
\end{equation}
where $(\alpha,\beta,\gamma,z_{12},z_{13},z_{23})$ are the Euler-Iwasawa coordinates attached to $k_0$.
\end{lemma}
\begin{remark}\label{independent}
(\ref{feid12}) is independent of $k_0$, which is also a consequence of $D_{12}$ being $G$-invariant. 
\end{remark}
{\noindent\bf Proof of Lemma \ref{eid12}.} According to Theorem \ref{HCiso}, 
we first compute the left-invariant differential $\widehat D_{12}$
on $SL_3(\mathbb R)$ associated with $D_{12}$ so that 
\begin{equation}
 \widehat D_{12}=\frac{1}{2}\left(R\left(E_{12}\right)R\left(E_{21}\right)+ R\left(E_{21}\right)R\left(E_{12}\right)\right).
\end{equation}

Consider the Euler-Iwasawa coordinates $(\alpha,\beta,\gamma$, $z_{12},z_{13},z_{23},\lambda_1,\lambda_2)$ attached to $k_0$ (see Definition \ref{eisl}). Then, by Lemmas \ref{r12}, \ref{r21} in Appendix \ref{lefti},
\begin{equation*}
\begin{split}
& \widehat D_{12}=\frac{1}{2}\left(\frac{\lambda_2}{\lambda_1}\left\{\sec\beta\sin\gamma\frac{\partial }{\partial \alpha}+\cos\gamma\frac{\partial }{\partial \beta}+\tan\beta\sin\gamma\frac{\partial}{\partial \gamma}+\left(z_{12}^2+1\right)\frac{\partial}{\partial z_{12}}+z_{23}\frac{\partial}{\partial z_{13}}-z_{13}\frac{\partial}{\partial z_{23}}\right\}\right.\\
&\,\,\,\,\,\,\,\,\,\,\left.+\lambda_2z_{12}\frac{\partial}{\partial \lambda_1}-\frac{\lambda_2^2}{\lambda_1}z_{12}\frac{\partial}{\partial \lambda_2}\right)\left(\frac{\lambda_1}{\lambda_2}\frac{\partial}{\partial z_{12}}\right)+\frac{1}{2}\left(\frac{\lambda_1}{\lambda_2}\frac{\partial}{\partial z_{12}}\right)\left(\frac{\lambda_2}{\lambda_1}\left\{\sec\beta\sin\gamma\frac{\partial }{\partial \alpha}+\cos\gamma\frac{\partial }{\partial \beta}\right.\right.\\
&\,\,\,\,\,\,\,\,\,\,\left.\left.+\tan\beta\sin\gamma\frac{\partial}{\partial \gamma}+\left(z_{12}^2+1\right)\frac{\partial}{\partial z_{12}}+z_{23}\frac{\partial}{\partial z_{13}}-z_{13}\frac{\partial}{\partial z_{23}}\right\}+\lambda_2z_{12}\frac{\partial}{\partial \lambda_1}-\frac{\lambda_2^2}{\lambda_1}z_{12}\frac{\partial}{\partial \lambda_2}\right)\\
&=\left(\sec\beta\sin\gamma\right)\frac{\partial^2 }{\partial \alpha\partial z_{12}}+\left(\cos\gamma\right)\frac{\partial^2}{\partial \beta\partial z_{12}} +\left(\tan\beta\sin\gamma\right) \frac{\partial^2}{\partial \gamma\partial z_{12}}+\left(z_{12}^2+1\right)\frac{\partial^2}{\partial z_{12}\partial z_{12}}\\
&\,\,\,\,\,\,\,\,\,\,\,\,\,+\left(z_{23}\right)\frac{\partial^2}{\partial z_{13}\partial z_{12}}+\left(-z_{13}\right)\frac{\partial^2}{\partial z_{23}\partial z_{12}}+2z_{12}\frac{\partial }{\partial z_{12}}\\
&\,\,\,\,\,\,\,\,\,\,\,\,\,+\left(\lambda_1z_{12}\right)\frac{\partial^2}{\partial z_{12}\partial\lambda_1}-\left(\lambda_2z_{12}\right)\frac{\partial^2}{\partial z_{12}\partial \lambda_2}+\left(\frac{\lambda_1}{2}\right)\frac{\partial}{\partial \lambda_1}-\left(\frac{\lambda_2}{2}\right) \frac{\partial}{\partial\lambda_2}\\
\end{split}
\end{equation*}
\normalsize

Let $f$ be an arbitary smooth function on $SL_3(\mathbb R)/A$ and $\widetilde f$ its lift to $SL_3(\mathbb R)$ via the natural projection $\pi:SL_3(\mathbb R)\rightarrow X$. Then,
\begin{equation*}
\begin{split}
&(D_{12}f)(\alpha,\beta,\gamma,z_{12},z_{13},z_{23})=\left(\widehat D_{12} \widetilde f\right)(\alpha,\beta,\gamma,z_{12},z_{13},z_{23},\lambda_1,\lambda_2)\\
&=\left(\sec\beta\sin\gamma\right)\frac{\partial^2 }{\partial \alpha\partial z_{12}}+\left(\cos\gamma\right)\frac{\partial^2}{\partial \beta\partial z_{12}} +\left(\tan\beta\sin\gamma\right) \frac{\partial^2}{\partial \gamma\partial z_{12}}+\left(z_{12}^2+1\right)\frac{\partial^2}{\partial z_{12}\partial z_{12}}\\
&\,\,\,\,\,\,\,\,\,\,\,\,\,+\left(z_{23}\right)\frac{\partial^2}{\partial z_{13}\partial z_{12}}+\left(-z_{13}\right)\frac{\partial^2}{\partial z_{23}\partial z_{12}}+2z_{12}\frac{\partial }{\partial z_{12}}.\\
\end{split}
\end{equation*}

We complete the proof of Lemma (\ref{eid12}). \,\,\,$\endpf$
\medskip

Next, we introduce a special class $\mathcal P$ of linear differential operators defined on $\underline{U}_{EI}\times \mathbb R^3$, in terms of the Euler-Iwasawa coordinates $(\alpha,\beta,\gamma,z_{12},z_{13}$, $z_{23})$. 
\begin{definition}\label{cp} The class $\mathcal P$ consists of 
\begin{equation*}
\sum\nolimits_{\substack{u_i\in U\left(\mathfrak{sl}_3(\mathbb R)\right)\\{\rm are\,\,monomials}}}\frac{P_i\left(\sin\alpha,\cos\alpha,\sin\beta,\cos\beta,\sin\gamma,\cos\gamma\right)}{\cos^{l_i}\beta}L_{u_i},
\end{equation*}
where $l_i\geq 0$, $P_i$ are polynomials in trigonometric functions,  $L_{u_i}$ are left derivatives on $X$.
\end{definition}
\begin{lemma}\label{composition}
If $D_1,D_2\in\mathcal P$, then $D_1\circ D_2\in\mathcal P$, where 
\begin{equation*}
(D_1\circ D_2)f:=D_1(D_2f),\,\,f\in C^{\infty}(X).
\end{equation*}
\end{lemma}
{\bf\noindent Proof of Lemma \ref{composition}.} Notice that $L_{X_{12}}$, $L_{X_{13}}$, $L_{X_{23}}$, $L_{X_{21}}$, $L_{X_{31}}$, $L_{X_{32}}$ constitute a set of generators of all left derivatives. By explicit formulas of  $L_{X_{ij}}$ in Lemmas \ref{x12}, \ref{x13}, \ref{x23}, \ref{x21}, \ref{x31}, \ref{x32} in Appendix \ref{LEFTD}, we can verify Lemma \ref{composition}.\,\,\,\,$\endpf$

\begin{lemma}\label{bounded}
Each differential operator $D\in\mathcal P$ can be written as 
\begin{equation*}
D=\sum\nolimits_{\substack{u_i\in U\left(\mathfrak{sl}_3(\mathbb R)\right)\\{\rm are\,\,monomials}}}r_iL_{u_i},
\end{equation*}
where $r_i$ are bounded functions on $\underline{U}_{EI}\times \mathbb R^3\subset SL_3(\mathbb R)/A$, and $L_{u_i}$ are left derivatives on $X$.
\end{lemma}
{\bf\noindent Proof of Lemma \ref{bounded}.}  It follows easily from Definition \ref{cp}.\,\,\,\,$\endpf$

\begin{lemma}\label{bk0}
The following differential operators are in $\mathcal P$.
\begin{equation*}
\frac{\partial}{\partial\alpha},\,\,\frac{\partial}{\partial\beta},\,\,\frac{\partial}{\partial\gamma},\,\,\frac{\partial }{\partial z_{13}},\,\,\frac{\partial }{\partial z_{23}},z_{12}\frac{\partial }{\partial z_{12}}-z_{23}\frac{\partial }{\partial z_{23}},\,\,\frac{\partial }{\partial z_{12}}+z_{23}\frac{\partial }{\partial z_{13}},\,\,z_{13}\frac{\partial }{\partial z_{13}}+z_{23}\frac{\partial }{\partial z_{23}}.
\end{equation*}
\end{lemma}
{\bf\noindent Proof of Lemma \ref{bk0}.} By Lemmas \ref{x21}, \ref{x31}, \ref{x32} in Appendix \ref{LEFTD}, we have \begin{equation}\label{p3}
\begin{split}
&\frac{\partial}{\partial\alpha}=\left(L_{X_{23}}-L_{X_{32}}\right),\,\,\,\frac{\partial}{\partial\beta}=\cos\alpha\left(L_{X_{12}}-L_{X_{21}}\right)-\sin\alpha\left(L_{X_{31}}-L_{X_{13}}\right),\\
&\frac{\partial}{\partial\gamma}=\cos\alpha\cos\beta\left(L_{X_{31}}-L_{X_{13}}\right)+\sin\alpha\cos\beta\left(L_{X_{12}}-L_{X_{21}}\right)-\sin\beta\left(L_{X_{23}}-L_{X_{32}}\right).
\end{split}	\end{equation}

Computation yields that	\begin{equation*}
\begin{split}
&\frac{1}{4}\left[\frac{\partial}{\partial\gamma},\left[\frac{\partial}{\partial\gamma},L_{X_{12}}\right]\right]+L_{X_{12}}-\left(\frac{1}{2}\sin\alpha\tan\beta\right)\frac{\partial }{\partial \alpha}-\left(\frac{1}{2}\cos\alpha\right)\frac{\partial }{\partial \beta}-\left(\frac{1}{2}\sin\alpha\sec\beta\right)\frac{\partial }{\partial \gamma}\\
&=-\left(\left(\frac{3}{4}\cos\alpha\sin 2\beta\right)z_{12}+\left(\frac{3}{4}\cos\alpha\cos 2\beta\cos\gamma-\frac{3}{4}\sin\alpha\sin\beta\sin\gamma\right)\right)\frac{\partial }{\partial z_{12}}\\
&\,\,\,\,\,\,\,-\left(\left(-\frac{3}{4}\sin\alpha\sin\beta\sin\gamma+\frac{3}{4}\cos\alpha\cos 2\beta\cos\gamma\right)z_{23}\right)\frac{\partial }{\partial z_{13}}\,\,\,\,\,\,\,\,\,\,\,\,\\
&\,\,\,\,\,\,\,-\left(\left(-\frac{3}{4}\cos\alpha\sin 2\beta\right)z_{23}+\left(\frac{3}{4}\sin\alpha\sin\beta\cos\gamma+\frac{3}{4}\cos\alpha\cos 2\beta\sin\gamma\right)\right)\frac{\partial }{\partial z_{23}}=:L^{\gamma\gamma}.\\
\end{split}
\end{equation*}
\normalsize
Then, by Lemma \ref{composition} $L^{\gamma\gamma}\in\mathcal P$. Similarly, 
\begin{equation*}
\begin{split}
&L^{\gamma\gamma\gamma\gamma}:=\left[\frac{\partial}{\partial\gamma},\left[\frac{\partial}{\partial\gamma},L^{\gamma\gamma}\right]\right]+L^{\gamma\gamma}=-\left(\frac{3}{4}\cos\alpha\sin 2\beta\right)\left(z_{12}\frac{\partial }{\partial z_{12}}-z_{23}\frac{\partial }{\partial z_{23}}\right)\in\mathcal P,\\
\end{split}
\end{equation*}
and hence
\begin{equation*}
\begin{split}
&\cos2\beta\left[\frac{\partial}{\partial\beta},\sin\alpha\left[\frac{\partial}{\partial\alpha},L^{\gamma\gamma\gamma\gamma}\right]-\cos\alpha\cdot L^{\gamma\gamma\gamma\gamma}\right]+2\sin2\beta\left(\sin\alpha\left[\frac{\partial}{\partial\alpha},L^{\gamma\gamma\gamma\gamma}\right]-\cos\alpha\cdot L^{\gamma\gamma\gamma\gamma}\right)\\
&\,\,\,\,\,\,\,\,\,\,\,=\frac{3}{2}\left(z_{12}\frac{\partial }{\partial z_{12}}-z_{23}\frac{\partial }{\partial z_{23}}\right)\in\mathcal P.\\
\end{split}
\end{equation*}
\normalsize
Compuation yields
\begin{equation*}
\begin{split}
&\widehat L:=L^{\gamma\gamma}+\left(\frac{3}{4}\cos\alpha\sin 2\beta\right)\left(z_{12}\frac{\partial}{\partial z_{12}}-z_{23}\frac{\partial }{\partial z_{23}}\right)\\
&\,\,\,\,\,\,=-\left(\frac{3}{4}\cos\alpha\cos 2\beta\cos\gamma-\frac{3}{4}\sin\alpha\sin\beta\sin\gamma\right)\left(\frac{\partial }{\partial z_{12}}+z_{23}\frac{\partial }{\partial z_{13}}\right)\\
&\,\,\,\,\,\,\,\,\,\,\,\,\,-\left(\frac{3}{4}\sin\alpha\sin\beta\cos\gamma+\frac{3}{4}\cos\alpha\cos 2\beta\sin\gamma\right)\frac{\partial }{\partial z_{23}},\\
&\widehat L^{\alpha}:=\sin\alpha\left[\frac{\partial}{\partial\alpha},\widehat L\right]-\cos\alpha\cdot\widehat L\\
&\,\,\,\,\,\,\,\,\,=\left(\frac{3}{4}\cos 2\beta\cos\gamma\right)\left(\frac{\partial }{\partial z_{12}}+z_{23}\frac{\partial }{\partial z_{13}}\right)+\left(\frac{3}{4}\cos 2\beta\sin\gamma\right)\frac{\partial }{\partial z_{23}},\\
\end{split}   
\end{equation*}
\begin{equation*}
\begin{split}
&\widehat L_1^{\alpha\gamma}:=-\sin\gamma\left[\frac{\partial}{\partial\gamma},\widehat L^{\alpha}\right]+\cos\gamma\cdot\widehat L^{\alpha}=\left(\frac{3}{4}\cos 2\beta\right)\left(\frac{\partial }{\partial z_{12}}+z_{23}\frac{\partial }{\partial z_{13}}\right),\\
&\widehat L_2^{\alpha\gamma}:=\cos\gamma\left[\frac{\partial}{\partial\gamma},\widehat L^{\alpha}\right]+\sin\gamma\cdot\widehat L^{\alpha}=\left(\frac{3}{4}\cos 2\beta\right)\frac{\partial }{\partial z_{23}},\\
\end{split}
\end{equation*}
and hence
\begin{equation*}
\begin{split}
&\widehat L_1^{\alpha\gamma\beta}:=-\sin2\beta\left[\frac{\partial}{\partial\beta},\widehat L_1^{\alpha\gamma}\right]+2\cos2\beta\cdot\widehat L_1^{\alpha\gamma}=\frac{3}{2}\left(\frac{\partial }{\partial z_{12}}+z_{23}\frac{\partial }{\partial z_{13}}\right)\in\mathcal P,\\
&\widehat L_2^{\alpha\gamma\beta}:=-\sin2\beta\left[\frac{\partial}{\partial\beta},\widehat L_2^{\alpha\gamma}\right]+2\cos2\beta\cdot\widehat L_2^{\alpha\gamma}=\frac{3}{2}\frac{\partial }{\partial z_{23}}\in\mathcal P.\\
\end{split}
\end{equation*}

Consider
\begin{equation*}
\begin{split}
&L_{X_{12}}+\left(\frac{1}{2}\cos\alpha\cos\beta\sin 2\gamma-\sin\alpha\tan\beta\sin^2\gamma-\frac{1}{2}\cos\alpha\sin\beta\tan\beta\sin 2\gamma\right)\frac{\partial }{\partial \alpha}\\
&\,\,\,\,\,\,+\left(-\frac{1}{2}\sin\alpha\sin\beta\sin 2\gamma-\cos\alpha\sin^2\beta\cos^2\gamma-\cos\alpha\cos^2\beta\sin^2\gamma\right)\frac{\partial }{\partial \beta}\\
&\,\,\,\,\,\,+\left(-\sin\alpha\sec\beta\sin^2\gamma-\frac{1}{2}\cos\alpha\tan\beta\sin 2\gamma+\frac{1}{4}\cos\alpha\sin 2\beta\sin 2\gamma \right)\frac{\partial }{\partial \gamma}\\
&\,\,\,\,\,\,+\left(\frac{1}{2}\sin\alpha\cos\beta\sin 2\gamma +\frac{1}{2}\cos\alpha\sin 2\beta\cos^2\gamma+\frac{1}{2}\cos\alpha\sin 2\beta\right)\left(z_{12}\frac{\partial }{\partial z_{12}}-z_{23}\frac{\partial }{\partial z_{23}}\right)\\
&\,\,\,\,\,\,+\left(-\sin\alpha\sin\beta\sin\gamma+\cos\alpha\cos 2\beta\cos\gamma\right)\left(\frac{\partial }{\partial z_{12}}+z_{23}\frac{\partial }{\partial z_{13}}\right)\\
&\,\,\,\,\,\,+\left(\sin\alpha\sin\beta\cos\gamma+\cos\alpha\cos 2\beta\sin\gamma\right)\frac{\partial }{\partial z_{23}}\\
&=-\left(\sin\alpha\cos\beta\sin 2\gamma+\frac{1}{2}\cos\alpha\sin 2\beta \cos 2\gamma\right)\left(z_{13}\frac{\partial }{\partial z_{13}}+z_{23}\frac{\partial }{\partial z_{23}}\right)\,\,\,\,\,\,\,\,\,\,\,\,\,\,\,\,\,\,\,\,\,\,\,\,\,\,\,\,\,\,\,\,\,\,\,\,\,\,\,\,\,\,\,\,\,\\
&\,\,\,\,\,\,\,-\left(-\sin\alpha\cos\beta\cos 2\gamma+\frac{1}{2}\cos\alpha \sin 2\beta\sin 2\gamma\right)\frac{\partial }{\partial z_{13}}.\\
\end{split}
\end{equation*}
By a similar argument, one can show $z_{13}\frac{\partial }{\partial z_{13}}+z_{23}\frac{\partial }{\partial z_{23}},\,\,\frac{\partial }{\partial z_{13}} \in\mathcal P$.
The proof of Lemma \ref{bk0} is complete.\,\,\,\,$\endpf$

\begin{lemma}\label{pp0} Denote the differentials $\frac{\partial }{\partial z_{13}},\frac{\partial }{\partial z_{23}},z_{12}\frac{\partial }{\partial z_{12}}-z_{23}\frac{\partial }{\partial z_{23}}$, $\frac{\partial }{\partial z_{12}}+z_{23}\frac{\partial }{\partial z_{13}}$, $z_{13}\frac{\partial }{\partial z_{13}}+z_{23}\frac{\partial }{\partial z_{23}}$ by $\delta_1,\delta_2,\delta_3,\delta_4,\delta_5$, respectively. Then, we can represent
\begin{equation*}
\delta_i=\sum\nolimits_{\substack{v_{i,j}\in U\left(\mathfrak{sl}_3(\mathbb R)\right)\\{\rm are\,\,monomials}}}p_{i,j}L_{v_{i,j}},\,\,\,1\leq i\leq 5,
\end{equation*}
where $p_{i,j}$ are bounded functions on $X$, and  $L_{v_{i,j}}$ are left derivatives on $X$.
\end{lemma}
{\bf\noindent Proof of Lemma \ref{pp0}.} We first note that the differential $\delta_i$ is globally defined on $X$. By Lemmas \ref{bk0}, \ref{bounded},   $\delta_i$ can be written on the open set $\underline{U}_{EI}\times \mathbb R^3\subset SL_3(\mathbb R)/A$ as
\begin{equation}\label{121}
\delta_i=\sum\nolimits_{\substack{v_{i,j}\in U\left(\mathfrak{sl}_3(\mathbb R)\right)\\{\rm are\,\,monomials}}}p_{i,j}L_{v_{i,j}}, \end{equation}
where $p_{i,j}$ are bounded functions on $\underline{U}_{EI}\times \mathbb R^3$, and  $L_{v_{i,j}}$ are left derivatives on $X$.
\medskip

{\noindent\bf Claim.} For $k_0\in SO_3(\mathbb R)$, $\delta_i$ can be written on the open set $\left(k_0\cdot\underline{U}_{EI}\right)\times \mathbb R^3$ as
\begin{equation*}
\delta_i=\sum\nolimits_{\substack{v_{i,j}^{k_0}\in U\left(\mathfrak{sl}_3(\mathbb R)\right)\\{\rm are\,\,monomials}}}p_{i,j}^{k_0}L_{v_{i,j}^{k_0}}, 
\end{equation*}
where $p_{i,j}^{k_0}$ are bounded functions on $\left(k_0\cdot\underline{U}_{EI}\right)\times \mathbb R^3$, and  $L_{v_{i,j}^{k_0}}$ are left derivatives on $X$.
\medskip

{\noindent\bf Proof of Claim.} Take the Euler-Iwasawa coordinates $(\alpha,\beta,\gamma,z_{12},z_{13},z_{23})$ attached to $k_0$. It is clear that the explicit formula of the differential $\delta$  in terms of the Euler-Iwasawa coordinates attached to $k_0$, coincides with that in terms of the Euler-Iwasawa coordinates attached to the identity of $SO_3(\mathbb R)$. 

Let Ad be the adjoint action of $SL_3(\mathbb R)$ on its lie algebra, that is, ${\rm Ad}(g)X=\left.\frac{d}{dt}\right|_{t=0}g\cdot\exp\left(tX\right)\cdot g^{-1}$.    
Then the adjoint action Ad acts on the universal enveloping algebra $U(\mathfrak{sl}_3(\mathbb R))$ by
\begin{equation*}
{\rm Ad}(g)(X_1\cdot X_2\cdot\cdots\cdot X_r)={\rm Ad}(g)X_1\cdot{\rm Ad}(g)X_2\cdot\cdots\cdot{\rm Ad}(g)X_r,
\end{equation*}
where $X_1,X_2,\cdots,X_r\in\mathfrak{sl}_3(\mathbb R)$. 
Since for $u\in U(\mathfrak{sl}_3(\mathbb{R}))$, the left derivatives $L_{{\rm Ad}(k_0)u}$ and $L_u$ are intertwined by the left $k_0$-translation, one sees that the explicit formula of $L_{{\rm Ad}(k_0)u}$ in terms of the Euler-Iwasawa coordinates attached to $k_0$, coincides with the explicit formula of  $L_{u}$ in terms of the Euler-Iwasawa coordinates attached to the identity of $SO_3(\mathbb R)$.
Now by (\ref{121}),   in terms of the Euler-Iwasawa coordinates attached to $k_0$, we can derive 
\begin{equation*}
\delta_i=\sum\nolimits_{\substack{v_{i,j}\in U\left(\mathfrak{sl}_3(\mathbb R)\right)\\{\rm are\,\,monomials}}}p_{i,j}L_{{\rm Ad}(k_0)v_{i,j}}, \end{equation*}
where $p_{i,j}$ and  $v_{i,j}$ are exactly the same as in (\ref{121}).

We complete the proof of Claim.\,\,\,$\endpf$
\medskip

Take points $k_1,\cdots,k_M\in SO_3(\mathbb R)$ such that $\left\{\left(k_l\cdot\underline{U}_{EI}\right)\right\}_{l=1}^M$ is an open cover of $SO_3(\mathbb R)$, and a partition of unity $\{\rho _{l}\}_{l=1}^M$
subordinate to $\left\{\left(k_l\cdot\underline{U}_{EI}\right)\right\}_{l=1}^M$. By Claim,  we have
\small\begin{equation*}
\begin{split}
&\delta_i=\sum_{l=1}^M\sum\nolimits_{\substack{v_{i,j}^{k_l}\in U\left(\mathfrak{sl}_3(\mathbb R)\right)\\{\rm are\,\,monomials}}}\left(\rho_l\cdot p_{i,j}^{k_l}\right)L_{v_{i,j}^{k_l}}.\\   
\end{split}
\end{equation*}
\normalsize
Noticing that the above   $\rho_l\cdot p_{i,j}^{k_l}$ are bounded on $X$, we  conclude Lemma \ref{pp0}.\,\,\,\, $\endpf$

\begin{lemma}\label{pp1} We can write
\begin{equation}\label{129}
\begin{split}
&D_{12}-\left(z_{12}^2+1\right)\frac{\partial^2}{\partial z_{12}\partial z_{12}}-\left(z_{23}\right)\frac{\partial^2}{\partial z_{13}\partial z_{12}}-\left(-z_{13}\right)\frac{\partial^2}{\partial z_{23}\partial z_{12}}-2z_{12}\frac{\partial }{\partial z_{12}}\\ 
&\,\,\,\,\,\,\,\,\,\,\,\,\,\,\,\,\,\,\,\,\,\,=\left(\sum_{\substack{w_j\in \mathfrak{so}_3(\mathbb R)}}q_jL_{w_j}\right) \frac{\partial}{\partial z_{12}}.\\
\end{split} \end{equation}
Here $q_j$ are bounded functions on $X$ which are independent of the variables $z_{12},z_{23},z_{13}$; $L_{w_j}$ are left derivatives on $X$ induced by the elements in $\mathfrak{so}_3(\mathbb R)$.
\end{lemma}

{\bf\noindent Proof of Lemma \ref{pp1}.}
The proof is similar to that of Lemma \ref{pp0}. By Lemma \ref{eid12}, in each coordinate chart $\left(\left(k_0\cdot\underline{U}_{EI}\right)\times\mathbb R^3,\Lambda_{k_0}\right)$,  
\begin{equation*}
	\begin{split}
		&D_{12}=\left(\sec\beta\sin\gamma\right)\frac{\partial^2 }{\partial \alpha\partial z_{12}}+\left(\cos\gamma\right)\frac{\partial^2}{\partial \beta\partial z_{12}} +\left(\tan\beta\sin\gamma\right) \frac{\partial^2}{\partial \gamma\partial z_{12}}\\
		&\,\,\,\,\,\,\,\,\,\,\,\,\,+\left(z_{12}^2+1\right)\frac{\partial^2}{\partial z_{12}\partial z_{12}}+\left(z_{23}\right)\frac{\partial^2}{\partial z_{13}\partial z_{12}}+\left(-z_{13}\right)\frac{\partial^2}{\partial z_{23}\partial z_{12}}+2z_{12}\frac{\partial }{\partial z_{12}},\\
	\end{split}
\end{equation*}
where $(\alpha,\beta,\gamma,z_{12},z_{13},z_{23})$ are the Euler-Iwasawa coordinates attached to $k_0\in SO_3(\mathbb R)$. By Lemma \ref{bounded} and (\ref{p3}) in Lemma \ref{bk0} , we can have (\ref{129}) on the open set $\underline{U}_{EI}\times\mathbb R^3\subset X$ such that $q_j$ are bounded functions on $\underline{U}_{EI}\times\mathbb R^3$ which are independent of the variables $z_{12},z_{23},z_{13}$, and $L_{w_j}$ are left derivatives on $X$ induced by the elements in $\mathfrak{so}_3(\mathbb R)$.

Take the Euler-Iwasawa coordinates attached to $k_0$.  Now by (\ref{129}),   in terms of the Euler-Iwasawa coordinates attached to $k_0$, 
\begin{equation*}
\begin{split}
&D_{12}-\left(z_{12}^2+1\right)\frac{\partial^2}{\partial z_{12}\partial z_{12}}-\left(z_{23}\right)\frac{\partial^2}{\partial z_{13}\partial z_{12}}-\left(-z_{13}\right)\frac{\partial^2}{\partial z_{23}\partial z_{12}}-2z_{12}\frac{\partial }{\partial z_{12}}\\ 
&\,\,\,\,\,\,\,\,\,\,\,\,\,\,\,\,\,\,\,\,\,\,=\left(\sum\nolimits_{\substack{w_j\in \mathfrak{so}_3(\mathbb R)}}q_jL_{{\rm Ad}(k_0)w_j}\right) \frac{\partial}{\partial z_{12}},\\
\end{split}
\end{equation*}
where $q_j$, $w_j$ are the same as in (\ref{129}). 
Take the partition of unity $\{\rho _{l}\}_{l=1}^M$
subordinate to $\left\{\left(k_l\cdot\underline{U}_{EI}\right)\right\}_{l=1}^M$ as in the proof of Lemma \ref{pp0}.
We can conclude Lemma \ref{pp1}.\,\,\,\, $\endpf$

\subsection{Proof
of the density and Theorem \ref{self}}\label{extra}

Assume that $X:=G/H$ is a homogeneous space with a $G$-invariant measure ${\rm d}\mu$.  Denote by $C^{\infty}(X)$ the space of smooth functions on $X$, and by $C^{\infty}_c(X)$ the space of functions in $C^{\infty}(X)$ with compact support. 

Assume $\Delta\in \mathbb D(X)$ is formally self-adjoint, that is, $(\Delta f,g)=(f,\Delta g)$ 
for all $f,g\in C_c^\infty (X)$, which also implies that the above equality also holds for all $f\in C_c^\infty(X)$ and $g$ in the space $ C_c^\infty(X)^*$ of distributions. It is clear that $\Delta$ is defined naturally in $\text{Dom}(\Delta):=\left\{f\in L^2(X): \Delta f\in L^2(X)\right\}$, and that the graph of the operator $\Delta: \text{Dom}(\Delta)\to L^2(X)$ is closed. The following two lemmas are standard. 
\begin{lemma}\label{dtos}
Suppose that $\Delta$ is a formally self-adjoint differential operator on $X$.
If the graph of $\Delta$ on the domain $C^{\infty}_c(X)$ is dense in the graph of $\Delta$ on the domain ${\rm Dom}(\Delta)$ with respect to the $L^2\times L^2$ norm, then $\Delta$ is symmetric on ${\rm Dom}(\Delta)$, that is,  
\begin{equation*}
(\Delta f,g)=(f,\Delta g) \,\,\,{\rm for\,\, all\,\,}f,g\in {\rm Dom}(\Delta).
\end{equation*}
\end{lemma}

\begin{lemma}\label{stos}
If the operator $\Delta$ is symmetric on ${\rm Dom}(\Delta)$, then $\Delta$ is self-adjoint 
on ${\rm Dom}(\Delta)$.
\end{lemma}

In the following, we shall establish the density of $C^{\infty}_c\left(SL_3(\mathbb R)/A\right)$ in ${\rm Dom}(D_{12})$ for the invariant differential operator $D_{12}$. 

Firstly, we construct a sequence of cutoff functions $\{\chi_n\}_{n=1}^{\infty}$ on $X$. Take $\xi(r)$ to be a smooth function on $[0,+\infty)$ such that $0\leq\xi(r)\leq 1$ for $r\in[0,+\infty)$, 
$\xi(r)\equiv 1$ for $0\leq r\leq 1$, and $\xi(r)\equiv 0$ for $r\geq 2$. For $n=1,2,\cdots$, define 
\begin{equation}\label{cut}
\chi_n(x):=\chi_n\left(k_x,z_{12},z_{13},z_{23}\right):=\xi\left(\frac{|z_{13}|^2+|z_{23}|^2}{n}\right)\cdot \xi\left(\frac{\ln\left(|z_{12}|^2+1\right)}{n}\right).
\end{equation}
\begin{lemma}\label{ncutoff}
The following properties hold for the cutoff functions $\{\chi_n\}_{n=1}^{\infty}$ defined by (\ref{cut}).
\begin{enumerate}[label={\rm(\alph*)}]
\item $\{\chi_n\}_{n=1}^{\infty}\subset C^{\infty}_c(X)$.
\item $0\leq\chi_n(x)\leq 1$, and $\lim_{n\to\infty}\chi_n(x)=1$, for each $x\in X$. 
\item For non-negative integers $i_1,i_2,i_3$ such that $i_1+i_2+i_3\geq 1$,
\begin{equation*}
 \lim_{n\to\infty}\frac{\partial^{i_1+i_2+i_3} }{\partial z_{12}^{i_1}\partial z_{13}^{i_2}\partial z_{23}^{i_3}}\chi_n(x)=0.
\end{equation*}
\item  For $n=1,2,\cdots$,  we have
\begin{equation*}
\frac{\partial \chi_n}{\partial \alpha}=\frac{\partial\chi_n}{\partial \beta}=\frac{\partial\chi_n}{\partial \gamma}\equiv 0.
\end{equation*}
\item  For $n=1,2,\cdots$,  we have
\begin{equation*}
\left(z_{23}\frac{\partial}{\partial z_{13}}-z_{13}\frac{\partial}{\partial z_{23}}\right)\chi_n\equiv 0.
\end{equation*}
\item There exists $\mathcal M>0$ such that for all $n=1,2,\cdots$,
\begin{equation*}
\sup_{x\in X}\left|(1
+|z_{12}|)(1+|z_{13}|+|z_{23}|)\cdot\frac{\partial \chi_n}{\partial z_{12}}(x)\right|+\sup_{x\in X}\left|\left(|z_{12}|^2+1\right)\cdot\frac{\partial^{2} \chi_n}{\partial z_{12}\partial z_{12}}(x)\right|\leq \mathcal M.
\end{equation*}
\end{enumerate}
\end{lemma}
{\bf\noindent Proof of Lemma \ref{ncutoff}.}
It is clear that (a), (b), (c), (d) follow from the definition in (\ref{cut}). (e) follows from the computation that
\begin{equation*}
z_{23}\frac{\partial \chi_n}{\partial z_{13}}-z_{13}\frac{\partial \chi_n}{\partial z_{23}}=\xi\left(\frac{\ln\left(|z_{12}|^2+1\right)}{n}\right)\cdot\xi^{\prime}\left(\frac{|z_{13}|^2+|z_{23}|^2}{n}\right)\cdot\left(\frac{z_{23}z_{13}}{n}-\frac{z_{13}z_{23}}{n}\right)\equiv 0. 
\end{equation*}

Since
\begin{equation*}
\begin{split}
&\frac{\partial\chi_n}{\partial z_{12}}\left(k_x,z_{12},z_{13},z_{23}\right)=\xi\left(\frac{|z_{13}|^2+|z_{23}|^2}{n}\right)\cdot \xi^{\prime}\left(\frac{\ln\left(|z_{12}|^2+1\right)}{n}\right)\cdot\frac{2z_{12}}{|z_{12}|^2+1}\cdot\frac{1}{n},\\
\end{split}
\end{equation*}
we have
\begin{equation*}
\begin{split}
&\sup_{x\in X}\left|(1
+|z_{12}|)(1+|z_{13}|+|z_{23}|)\cdot\frac{\partial \chi_n}{\partial z_{12}}(x)\right|\\
&=\sup_{x\in X}\left(\left|\xi\left(\frac{|z_{13}|^2+|z_{23}|^2}{n}\right) \xi^{\prime}\left(\frac{\ln\left(|z_{12}|^2+1\right)}{n}\right)\right|\cdot\frac{2|z_{12}|(1
+|z_{12}|)}{|z_{12}|^2+1}\cdot\frac{1+|z_{13}|+|z_{23}|}{n}\right)\\
&\leq 4\sup_{r\in[0,+\infty)}\left|\xi^{\prime}(r)\right|\cdot\sup_{x\in X}\left(\left|\xi\left(\frac{|z_{13}|^2+|z_{23}|^2}{n}\right)\right|\cdot\frac{1+|z_{13}|+|z_{23}|}{n}\right).\\
\end{split}
\end{equation*}
Since when $1+|z_{13}|+|z_{23}|\geq4n$, 
$\xi\left(\frac{|z_{13}|^2+|z_{23}|^2}{n}\right)=0$, we conclude that
\begin{equation*}
\begin{split}
&\sup_{x\in X}\left|(1
+|z_{12}|)(1+|z_{13}|+|z_{23}|)\cdot\frac{\partial \chi_n}{\partial z_{12}}(x)\right|\leq 16\sup_{r\in[0,+\infty)}\left|\xi^{\prime}(r)\right|.\\
\end{split}
\end{equation*}

Similarly, we have 
\begin{equation*}
\begin{split}
&\frac{\partial^2\chi_n}{\partial z_{12}\partial z_{12}}\left(k_x,z_{12},z_{13},z_{23}\right)=\xi\left(\frac{|z_{13}|^2+|z_{23}|^2}{n}\right)\cdot\\
&\cdot\left(\xi^{\prime\prime}\left(\frac{\ln\left(|z_{12}|^2+1\right)}{n}\right)\cdot\frac{4z_{12}^2}{(|z_{12}|^2+1)^2}\cdot\frac{1}{n^2}+\xi^{\prime}\left(\frac{\ln\left(|z_{12}|^2+1\right)}{n}\right)\cdot\frac{2(1-z_{12}^2)}{(|z_{12}|^2+1)^2}\cdot\frac{1}{n}\right),\\
\end{split}
\end{equation*}
and hence
\begin{equation*}
\begin{split}
&\sup_{x\in X}\left|\left(|z_{12}|^2+1\right)\cdot\frac{\partial^{2} \chi_n}{\partial z_{12}^{2}}(x)\right|\leq 4\sup_{r\in[0,+\infty)}\left(\left|\xi^{\prime}(r)\right|+\left|\xi^{\prime\prime}(r)\right|\right),\\ 
\end{split}
\end{equation*}
We establish (f) and thus complete the proof of Lemma \ref{ncutoff}.\,\,\,\,$\endpf$
\medskip

Let $\{\psi_m\}_{m=1}^{\infty}$ be a smooth approximate identity at the identity element $\rm Id$  of $SL_3(\mathbb R)$  satisfying the following properties. $\psi_m$ are smooth non-negative functions on $SL_3(\mathbb R)$; $\psi_m$, are compactly supported in the open set $U_m\subset SL_3(\mathbb R)$, where $\cap_{m=1}^{\infty}U_m=\{\rm Id\}$; $\int _{SL_3(\mathbb R)}\psi_m (g)\,\mathrm {d} g=1$, where $\mathrm {d} g$ is a fixed Haar measure on $SL_3(\mathbb R)$.

\begin{lemma}\label{moll} The following holds for each $f\in L^2(X)$.
\begin{enumerate}[label={\rm (\alph*)}]
    \item $\psi_m*f$ are smooth functions on $X$ for $m=1,2,\cdots$, where
    \begin{equation}\label{conv}
(\psi_m*f)(x):=\int_{SL_3(\mathbb R)}\psi_m(g)\cdot f\left(g^{-1}x\right){\rm d}g. 
\end{equation} 
    
    \item $\psi_m*f\in L^2(X)$ for $m=1,2,\cdots$, and $\psi_m*f\rightarrow f$ in $L^2(X)$ as $m\rightarrow\infty$.
    
    \item  For each $u\in U\left(\mathfrak{sl}_3(\mathbb R)\right)$, $L_u(\psi_m*f)\in L^2(X)$.  
\end{enumerate}
Assume that $\Delta$ is an $SL_3(\mathbb R)$-invariant differential operator on $X$, and that $f,\Delta f\in L^2(X)$, then
\begin{enumerate}[label={\rm (\alph*)}]
\setcounter{enumi}{3}
    \item  $\Delta(\psi_m*f)\in L^2(X)$ for $m=1,2,\cdots$, and $\Delta(\psi_m*f)\rightarrow\Delta f$ in $L^2(X)$ as $m\rightarrow\infty$.
\end{enumerate}
\end{lemma}
{\noindent\bf Proof of Lemma \ref{moll}.}
See Appendix \ref{REG}.\,\,\,$\endpf$ 

\medskip

Now, we proceed to prove

\begin{theorem}\label{denass}
$D_{12}$ is an $SL_3(\mathbb R)$-invariant, formally self-adjoint differential operator on $X$. The graph of $D_{12}$ on the domain $C^{\infty}_c(X)$ is dense in the graph of $D_{12}$ on the domain ${\rm Dom}(D_{12})$, with respect to the graph $L^2\times L^2$ norm.
\end{theorem}
{\bf\noindent Proof of Theorem \ref{denass}.} Firstly, we show that $D_{12}$ is formally self-adjoint. It suffices to prove in each coordinate chart that $D_{12}$ coincides with $D_{12}^*$. This can be done by a direct computation, thanks to the invariant measure given by (\ref{hm}) and the explicit formula for $D_{12}$ given in Lemma \ref{eid12} (or see
Remark \ref{cfree} for a coordinate-free approach).

For each $f\in {\rm Dom}(D_{12})$, consider the functions $\chi_n\cdot(\psi_m*f)$, $m,n=1,2,\cdots$, where $\chi_n$ and $(\psi_m*f)$ are defined by (\ref{cut}) and (\ref{conv}), respectively. It is clear that $\chi_n\cdot(\psi_m*f)\in C^{\infty}_c(X)$ by (a) in Lemma \ref{ncutoff} and (a) in Lemma \ref{moll}. Further, by (b) and (d) in Lemma \ref{moll}, to prove the density of $C^{\infty}_c(X)$ in ${\rm Dom}(D_{12})$, it suffices to prove that for fixed $m=1,2,\cdots$,
\begin{equation}\label{pl2}
\begin{split}
&\lim_{n\to\infty}\|\chi_n\cdot(\psi_m*f)-\psi_m*f\|_{L^2(X)}=0,\\
\end{split}   
\end{equation}
and
\begin{equation}\label{pdl2}
\begin{split}
&\lim_{n\to\infty}\|D_{12}\left(\chi_n\cdot(\psi_m*f)\right)-D_{12}(\psi_m*f)\|_{L^2(X)}=0.\\
\end{split}   
\end{equation}
(\ref{pl2}) follows easily from (b) in Lemma \ref{ncutoff}, and Lebesgue's dominated convergence theorem. In the following, we shall establish (\ref{pdl2})  for $m=1,2,\cdots$. 

By Lemma \ref{pp1}, computation yields that
\begin{equation*}
\begin{split}
&D_{12}\left(\chi_n\cdot\left(\psi_m*f\right)\right)=\chi_n\cdot D_{12}\left(\psi_m*f\right)+\left(\sum\nolimits_{\substack{w_j\in \mathfrak{so}_3(\mathbb R)}}q_jL_{w_j}\left(\frac{\partial\chi_n}{\partial z_{12}}\right)\right)\cdot \left(\psi_m*f\right)\\
&+\left(\sum\nolimits_{\substack{w_j\in \mathfrak{so}_3(\mathbb R)}}q_jL_{w_j}\chi_n\right)\cdot \frac{\partial\left(\psi_m*f\right)}{\partial z_{12}}+\frac{\partial\chi_n}{\partial z_{12}}\cdot \left(\sum\nolimits_{\substack{w_j\in \mathfrak{so}_3(\mathbb R)}}q_jL_{w_j}\left(\psi_m*f\right)\right)\\
&+\left(z_{12}^2+1\right)\frac{\partial^2\chi_n}{\partial z_{12}\partial z_{12}}\cdot\left(\psi_m*f\right)+2\left(z_{12}^2+1\right)\frac{\partial\chi_n}{\partial z_{12}}\cdot\frac{\partial\left(\psi_m*f\right)}{\partial z_{12}}\\
&+\left(z_{23}\frac{\partial}{\partial z_{13}}-z_{13}\frac{\partial}{\partial z_{23}}\right)\left(\frac{\partial\chi_n}{\partial z_{12}}\right)\cdot\left(\psi_m*f\right)+\left(z_{23}\frac{\partial \chi_n}{\partial z_{13}}-z_{13}\frac{\partial\chi_n}{\partial z_{23}}\right)\cdot\frac{\partial\left(\psi_m*f\right)}{\partial z_{12}}\\
&+\frac{\partial\chi_n}{\partial z_{12}}\cdot\left(z_{23}\frac{\partial}{\partial z_{13}}-z_{13}\frac{\partial}{\partial z_{23}}\right)\left(\psi_m*f\right)+2z_{12}\frac{\partial\chi_n }{\partial z_{12}}\cdot\left(\psi_m*f\right).\\
\end{split}    
\end{equation*}
Since $L_{w_j}$,  $w_j\in\mathfrak{so}_3(\mathbb R)$, are linear combinations of $\frac{\partial}{\partial\alpha}$, $\frac{\partial}{\partial\beta}$, $\frac{\partial}{\partial\gamma}$, by (d) in Lemma \ref{ncutoff},
\begin{equation*}
\left(\sum\nolimits_{\substack{w_j\in \mathfrak{so}_3(\mathbb R)}}q_jL_{w_j}\left(\frac{\partial\chi_n}{\partial z_{12}}\right)\right)\cdot \left(\psi_m*f\right)=\left(\sum\nolimits_{\substack{w_j\in \mathfrak{so}_3(\mathbb R)}}q_jL_{w_j}\chi_n\right)\cdot \frac{\partial\left(\psi_m*f\right)}{\partial z_{12}}\equiv0.
\end{equation*}
Similarly, by (e) in Lemma \ref{ncutoff},
\begin{equation*}
\left(z_{23}\frac{\partial}{\partial z_{13}}-z_{13}\frac{\partial}{\partial z_{23}}\right)\left(\frac{\partial\chi_n}{\partial z_{12}}\right)\cdot\left(\psi_m*f\right)=\left(z_{23}\frac{\partial \chi_n}{\partial z_{13}}-z_{13}\frac{\partial\chi_n}{\partial z_{23}}\right)\cdot\frac{\partial\left(\psi_m*f\right)}{\partial z_{12}}\equiv0.
\end{equation*}
Therefore,
\begin{equation*}
\begin{split}
&D_{12}\left(\chi_n\left(\psi_m*f\right)\right)=\chi_nD_{12}\left(\psi_m*f\right)+\frac{\partial\chi_n}{\partial z_{12}} \left(\sum\nolimits_{\substack{w_j\in \mathfrak{so}_3(\mathbb R)}}q_jL_{w_j}\left(\psi_m*f\right)\right)\\
&+\left(z_{12}^2+1\right)\frac{\partial^2\chi_n}{\partial z_{12}\partial z_{12}}\left(\psi_m*f\right)+2\left(z_{12}^2+1\right)\frac{\partial\chi_n}{\partial z_{12}}\frac{\partial\left(\psi_m*f\right)}{\partial z_{12}}\\
&+\frac{\partial\chi_n}{\partial z_{12}}\cdot\left(z_{23}\frac{\partial}{\partial z_{13}}-z_{13}\frac{\partial}{\partial z_{23}}\right)\left(\psi_m*f\right)+2z_{12}\frac{\partial\chi_n }{\partial z_{12}}\cdot\left(\psi_m*f\right)\\
&=:A_n+B_n+C_n+D_n+E_n+F_n.\\
\end{split}    
\end{equation*}

In what follows, we shall prove that the functions $B_n$, $C_n$, $D_n$, $E_n$, $F_n$ converge to $0$ and $A_n$ converges to $D_{12}(\psi_m*f)$  in $L^2(X)$, when $n$ approaches $\infty$. Thanks to (f) in Lemma \ref{ncutoff} and the fact that $q_j$ are bounded functions on $X$, we conclude that
\begin{equation*}
\begin{split}
|B_n(x)|\leq&\left(\mathcal M\sup\nolimits_{\substack{w_j\in \mathfrak{so}_3(\mathbb R)\\x\in X}}\left|q_j(x)\right|\right)\sum\nolimits_{w_j\in \mathfrak{so}_3(\mathbb R)} \left|L_{w_j}\left(\psi_m*f\right)(x)\right|,\\\
\end{split}
\end{equation*}
where the right-hand side is in $L^2(X)$ by (c) in Lemma \ref{moll}. Then by (c) in Lemma \ref{ncutoff}, we  conclude by Lebesgue's dominated convergence theorem that $\lim_{n\to\infty}\|B_n\|_{L^2(X)}=0$. Similarly, by (f) in Lemma \ref{ncutoff}, we have pointwisely
\begin{equation*}
\begin{split}
|C_n(x)|+|F_n(x)|\leq 2\mathcal M\cdot \left|\left(\psi_m*f\right)(x)\right|.\\
\end{split}
\end{equation*}
Then by (b) in Lemma \ref{moll} and (c) in Lemma \ref{ncutoff}, we have
$\lim_{n\to\infty}\|C_n\|_{L^2(X)}=0$, $\lim_{n\to\infty}\|F_n\|_{L^2(X)}=0$.

By Lemma \ref{pp0}, we have
\small\begin{equation*}
\begin{split}
&\left(z_{12}^2+1\right)\frac{\partial\chi_n}{\partial z_{12}}\frac{\partial\left(\psi_m*f\right)}{\partial z_{12}}=\frac{\partial\chi_n}{\partial z_{12}}\cdot\left(\frac{\partial}{\partial z_{12}}+z_{23}\frac{\partial}{\partial z_{13}}\right)\left(\psi_m*f\right)-z_{23}\frac{\partial\chi_n}{\partial z_{12}}\frac{\partial\left(\psi_m*f\right)}{\partial z_{13}}\\  
&+z_{12}\frac{\partial\chi_n}{\partial z_{12}}\cdot\left(z_{12}\frac{\partial}{\partial z_{12}}-z_{23}\frac{\partial}{\partial z_{23}}\right)\left(\psi_m*f\right)+z_{12}z_{23}\frac{\partial\chi_n}{\partial z_{12}}\frac{\partial\left(\psi_m*f\right)}{\partial z_{23}}\\
&=\frac{\partial\chi_n}{\partial z_{12}}\left(\sum_{\substack{v_{4,j}\in U\left(\mathfrak{sl}_3(\mathbb R)\right)\\{\rm are\,\,monomials}}}p_{4,j}L_{v_{4,j}}\left(\psi_m*f\right)\right)-z_{23}\frac{\partial\chi_n}{\partial z_{12}}\left(\sum_{\substack{v_{1,j}\in U\left(\mathfrak{sl}_3(\mathbb R)\right)\\{\rm are\,\,monomials}}}p_{1,j}L_{v_{1,j}}\left(\psi_m*f\right)\right)\\
&+z_{12}\frac{\partial\chi_n}{\partial z_{12}}\left(\sum_{\substack{v_{3,j}\in U\left(\mathfrak{sl}_3(\mathbb R)\right)\\{\rm are\,\,monomials}}}p_{3,j}L_{v_{3,j}}\left(\psi_m*f\right)\right)+z_{12}z_{23}\frac{\partial\chi_n}{\partial z_{12}}\left(\sum_{\substack{v_{2,j}\in U\left(\mathfrak{sl}_3(\mathbb R)\right)\\{\rm are\,\,monomials}}}p_{2,j}L_{v_{2,j}}\left(\psi_m*f\right)\right).\\
\end{split}   
\end{equation*}
\normalsize
By (f) in Lemma \ref{ncutoff} and the fact that $p_{1,j},\cdots,p_{4,j}$ are bounded functions on $X$, 
\begin{equation*}
\begin{split}
|D_n(x)|\leq&\left(\mathcal M\sup\nolimits_{\substack{1\leq i\leq 4\\x\in X}}\left|p_{i,j}(x)\right|\right)\sum\nolimits_{i=1}^4\sum\nolimits_{\substack{v_{i,j}\in U\left(\mathfrak{sl}_3(\mathbb R)\right)\\{\rm are\,\,monomials}}}|L_{v_{i,j}}\left(\psi_m*f\right)(x)|.\\
\end{split}
\end{equation*}
By (c) in Lemma \ref{ncutoff}, we have $\lim_{n\to\infty}\|D_n\|_{L^2(X)}=0$. Similarly, by Lemma \ref{pp0},
\begin{equation*}
\begin{split}
&\frac{\partial\chi_n}{\partial z_{12}}\cdot\left(z_{23}\frac{\partial}{\partial z_{13}}-z_{13}\frac{\partial}{\partial z_{23}}\right)\left(\psi_m*f\right)=z_{23}\frac{\partial\chi_n}{\partial z_{12}}\cdot\left(\sum\nolimits_{\substack{v_{1,j}\in U\left(\mathfrak{sl}_3(\mathbb R)\right)\\{\rm are\,\,monomials}}}p_{4,j}L_{v_{4,j}}\left(\psi_m*f\right)\right)\\\\ 
&\,\,\,\,\,\,\,\,\,\,\,-z_{13}\frac{\partial\chi_n}{\partial z_{12}}\cdot\left(\sum\nolimits_{\substack{v_{2,j}\in U\left(\mathfrak{sl}_3(\mathbb R)\right)\\{\rm are\,\,monomials}}}p_{2,j}L_{v_{2,j}}\left(\psi_m*f\right)\right).\\
\end{split}
\end{equation*}
By (c), (f) in Lemma \ref{ncutoff}, and (c) in Lemma \ref{moll}, we have
$\lim_{n\to\infty}\|E_n\|_{L^2(X)}=0$.

Finally, we consider $A_n$. It is clear that
\begin{equation*}
\begin{split}
&|A_n(x)-\left(D_{12}\left(\psi_m*f\right)\right)(x)|=|\left(\chi_n-1\right)\cdot\left(D_{12}\left(\psi_m*f\right)\right)(x)|\leq|\left(D_{12}\left(\psi_m*f\right)\right)(x)|.\\    
\end{split}    
\end{equation*}
By (d) in Lemma \ref{moll} and (b) in Lemma \ref{ncutoff},
$\lim_{n\to\infty}\|A_n-D_{12}\left(\psi_m*f\right)\|_{L^2(X)}=0$. 

We complete the proof of Theorem \ref{denass}.\,\,\,\,$\endpf$
\medskip

{\noindent\bf Proof of Theorem \ref{self}.} By Theorem \ref{denass}, the graph of $D_{12}$ on the domain $C^{\infty}_c(X)$ is dense in the graph of $D_{12}$ on the domain $\text{Dom}(D_{12})$, with respect to the graph $L^2\times L^2$ norm. Then $D_{12}$ is essentially self-adjoint by Lemmas \ref{dtos}, \ref{stos}.

Notice that the conjugations of the normalizer of $A$ in $SL_{3}(\mathbb R)$ (or equivalently the  permutations of the indices $1,2,3$) are isometries of $L^2(X)$, and such isometries transform $D_{12}$ to $D_{13}$ and $D_{23}$. Then, $D_{13}$ and $D_{23}$ are essentially self-adjoint. 

We complete the proof of Theorem \ref{self}. \,\,\,\, $\endpf$

\begin{appendices}

\section{Computations in the Euler-Iwasawa Coordinates}
\subsection{Explicit formulas for the generators of the left derivatives on \texorpdfstring{$SL_3(\mathbb R)/A$}{dd}}\label{LEFTD}

\begin{lemma}\label{x12} In the coordinate chart $\left(\left(k_0\cdot\underline{U}_{EI}\right)\times\mathbb R^3,\Lambda_{k_0}\right)$, where $k_0$ is the identity of $SO_3(\mathbb R)$, 
\begin{equation*}
\begin{split}
&-L_{X_{12}}=\left(\frac{1}{2}\cos\alpha\cos\beta\sin 2\gamma-\sin\alpha\tan\beta\sin^2\gamma-\frac{1}{2}\cos\alpha\sin\beta\tan\beta\sin 2\gamma\right)\frac{\partial }{\partial \alpha}\\
&+\left(-\frac{1}{2}\sin\alpha\sin\beta\sin 2\gamma-\cos\alpha\sin^2\beta\cos^2\gamma-\cos\alpha\cos^2\beta\sin^2\gamma\right)\frac{\partial }{\partial \beta}\\
&+\left(-\sin\alpha\sec\beta\sin^2\gamma-\frac{1}{2}\cos\alpha\tan\beta\sin 2\gamma+\frac{1}{4}\cos\alpha\sin 2\beta\sin 2\gamma \right)\frac{\partial }{\partial \gamma}\\
&+\left(\left(\frac{1}{2}\sin\alpha\cos\beta\sin 2\gamma +\frac{1}{2}\cos\alpha\sin 2\beta\cos^2\gamma+\frac{1}{2}\cos\alpha\sin 2\beta\right)z_{12}\right.\\
&\,\,\,\,\,\,\,\,\,\,\,\,+\left(-\sin\alpha\sin\beta\sin\gamma+\cos\alpha\cos 2\beta\cos\gamma\right)\biggr)\frac{\partial }{\partial z_{12}}\\
&+\left(\left(\sin\alpha\cos\beta\sin 2\gamma+\frac{1}{2}\cos\alpha\sin 2\beta \cos 2\gamma\right)z_{13}+\left(-\sin\alpha\sin\beta\sin\gamma\right.\right.\\
&\,\,\,\,\,\,\,\,\,\,\,\,\left.\left.+\cos\alpha\cos 2\beta\cos\gamma\right)z_{23}+\left(-\sin\alpha\cos\beta\cos 2\gamma+\frac{1}{2}\cos\alpha \sin 2\beta\sin 2\gamma\right)\right)\frac{\partial }{\partial z_{13}}\,\,\,\,\,\,\,\,\,\,\,\,\\
&+\left(\left(\frac{1}{2}\sin\alpha\cos\beta\sin 2\gamma-\frac{1}{2}\cos\alpha\sin 2\beta-\frac{1}{2}\cos\alpha\sin 2\beta\sin^2\gamma\right)z_{23}\right.\\
&\,\,\,\,\,\,\,\,\,\,\,\,+\left(\sin\alpha\sin\beta\cos\gamma+\cos\alpha\cos 2\beta\sin\gamma\right)\biggr)\frac{\partial }{\partial z_{23}}.\\
	\end{split}
\end{equation*}
\end{lemma}
{\noindent\bf Proof of Lemma \ref{x12}.} According to the definition,
\begin{equation*}
\begin{split}
&(-L_{X_{12}}f)(\alpha,\beta,\gamma,z_{12},z_{13},z_{23})=\left.\frac{d}{dt}\left(f\left(\widetilde \alpha(t),\widetilde\beta(t),\widetilde\gamma(t),\widetilde z_{12}(t),\widetilde z_{13}(t),\widetilde z_{23}(t)\right)\right)\right|_{t=0},\\
\end{split}
\end{equation*}
where $\left(\widetilde \alpha(t),\widetilde\beta(t),\widetilde\gamma(t),\widetilde z_{12}(t),\widetilde z_{13}(t),\widetilde z_{23}(t)\right)$ are the Euler-Iwasawa coordinates of
\begin{equation*}
\begin{split}
\left(\begin{matrix}1&t&0\\		0&1&0\\	0&0&1\\\end{matrix}\right)&\left(\begin{matrix}\cos{\beta}\cos{\gamma}&-\sin{\beta}&\cos{\beta}\sin{\gamma}\\\sin{\alpha}\sin{\gamma}+\cos{\alpha}\cos{\gamma}\sin{\beta}&\cos{\alpha}\cos{\beta}&\cos{\alpha}\sin{\beta}\sin{\gamma}-\cos{\gamma}\sin{\alpha}\\\cos{\gamma}\sin{\alpha}\sin{\beta}-\cos{\alpha}\sin{\gamma}&\cos{\beta}\sin{\alpha}&\cos{\alpha}\cos{\gamma}+\sin{\alpha}\sin{\beta}\sin{\gamma}\end{matrix}\right)\left(\begin{matrix}1&z_{12}&z_{13}\\		0&1&z_{23}\\
	0&0&1\\\end{matrix}\right)\\
&=:\left(\begin{matrix}1&t&0\\		0&1&0\\	0&0&1\\\end{matrix}\right)\left(\begin{matrix}k_{11}&k_{12}&k_{13}\\		k_{21}&k_{22}&k_{23}\\		k_{31}&k_{32}&k_{33}\end{matrix}\right)\left(\begin{matrix}1&z_{12}&z_{13}\\		0&1&z_{23}\\
	0&0&1\\\end{matrix}\right).\\
\end{split}
\end{equation*}
\normalsize
Applying the Gram-Schmidt process, we can obtain that,  up to order $1$ in $t$, 
\begin{equation*}
\begin{split}
\left(\begin{matrix}1&t&0\\		0&1&0\\	0&0&1\\\end{matrix}\right)\left(\begin{matrix}k_{11}&k_{12}&k_{13}\\		k_{21}&k_{22}&k_{23}\\		k_{31}&k_{32}&k_{33}\end{matrix}\right)\left(\begin{matrix}1&z_{12}&z_{13}\\		0&1&z_{23}\\
	0&0&1\\\end{matrix}\right)=\left(\begin{matrix}\widetilde k_{11}&\widetilde k_{12}&\widetilde k_{13}\\		\widetilde k_{21}&\widetilde k_{22}&\widetilde k_{23}\\		\widetilde k_{31}&\widetilde k_{32}&\widetilde k_{33}\end{matrix}\right)\left(\begin{matrix}1&\widetilde z_{12}&\widetilde z_{13}\\		0&1&\widetilde z_{23}\\
	0&0&1\\\end{matrix}\right)\left(\begin{matrix}\Delta_1&0&0\\	0&\Delta_2&0\\
	0&0&\Delta_3\end{matrix}\right),\\
\end{split}
\end{equation*}
\normalsize
where
\begin{equation*}
\begin{split}
&\Delta_1=1+t\left(k_{11}k_{21}\right),\,\,\Delta_2=1+t\left(k_{12}k_{22}\right),\,\,\Delta_3=1+t\left(k_{13}k_{23}\right),\\
&\widetilde k_{11}=k_{11}+t\left(k_{21}-k_{11}^2k_{21}\right),\,\,\widetilde k_{21}=k_{21}+t\left(-k_{11}k_{21}^2\right),\,\,\widetilde k_{31}=k_{31}+t\left(-k_{11}k_{21}k_{31}\right),\\
&\widetilde k_{12}=k_{12}+t\left(k_{13}^2k_{22}-k_{11}k_{12}k_{21}\right),\,\,\\
&\widetilde k_{22}=k_{22}+t\left(-k_{12}k_{21}^2+k_{13}k_{23}k_{22}\right),\\		&\widetilde k_{32}=k_{32}+t\left(k_{12}k_{23}k_{33}-k_{11}k_{22}k_{31}\right),\,\,\\
&\widetilde k_{13}=k_{13}+t\left(k_{13}^2k_{23}\right),\,\,\widetilde k_{23}=k_{23}+t\left(-k_{13} +k_{13}k_{23}^2\right),\,\,\widetilde k_{33}=k_{33}+t\left(k_{13}k_{23}k_{33}\right),\\	
&\widetilde z_{12}=\left(1+t\left(k_{11}k_{21}-k_{12}k_{22}\right)\right)z_{12}+t\left(k_{11}k_{22}+k_{12}k_{21}\right),\\
&\widetilde z_{13}=\left(1+t\left(k_{11}k_{21}-k_{13}k_{23}\right)\right)z_{13}+t\left(k_{11}k_{22}+k_{12}k_{21}\right)z_{23}+t\left(k_{11}k_{23}+k_{13}k_{21}\right),\\
&\widetilde z_{23}=\left(1+t\left(k_{12}k_{22}-k_{13}k_{23}\right)\right)z_{23}+t\left(k_{12}k_{23}+k_{13}k_{22}\right).\\
\end{split}
\end{equation*}
Then for smooth functions $f$ on $X$, we have by (\ref{inverse}) that
\begin{equation*}
\begin{split}
&(-L_{X_{12}}f)(\alpha,\beta,\gamma,z_{12},z_{13},z_{23})=\frac{\partial f}{\partial \alpha}\frac{d\widetilde\alpha}{dt}(0)+\frac{\partial f}{\partial \beta}\frac{d\widetilde\beta}{dt}(0)+\frac{\partial f}{\partial \gamma}\frac{d\widetilde\gamma}{dt}(0)+\frac{\partial f}{\partial z_{12}}\frac{d\widetilde z_{12}}{dt}(0)\\	
&\,\,\,\,\,\,\,\,\,\,\,\,\,\,\,\,\,\,\,\,\,\,\,\,+\frac{\partial f}{\partial z_{13}}\frac{d\widetilde z_{13}}{dt}(0)+\frac{\partial f}{\partial z_{23}}\frac{d\widetilde z_{23}}{dt}(0)\\
&=\left(\frac{(k_{21}k_{32}-k_{22}k_{31})(k_{12}k_{21}+k_{11}k_{22})}{k_{22}^2+k_{32}^2}\right)\frac{\partial f}{\partial \alpha}+\left(\frac{k_{11}k_{12}k_{21}-k_{13}^2k_{22}}{\sqrt{1-k_{12}^2}}\right)\frac{\partial f}{\partial \beta}\\
&\,\,\,\,\,\,\,\,\,\,\,\,\,+\left(\frac{-k_{13}k_{21}-k_{11}k_{12}k_{13}k_{22}}{k_{11}^2+k_{13}^2}\right)\frac{\partial f}{\partial \gamma}\\
&\,\,\,\,\,\,\,\,\,\,\,\,\,+\left(\left(k_{11}k_{21}-k_{12}k_{22}\right)z_{12}+\left(k_{11}k_{22}+k_{12}k_{21}\right)\right)\frac{\partial f}{\partial z_{12}}\\
&\,\,\,\,\,\,\,\,\,\,\,\,\,+\left(\left(k_{11}k_{21}-k_{13}k_{23}\right)z_{13}+\left(k_{11}k_{22}+k_{12}k_{21}\right)z_{23}+\left(k_{11}k_{23}+k_{13}k_{21}\right)\right)\frac{\partial f}{\partial z_{13}}\\
&\,\,\,\,\,\,\,\,\,\,\,\,\,+\left(\left(k_{12}k_{22}-k_{13}k_{23}\right)z_{23}+\left(k_{12}k_{23}+k_{13}k_{22}\right)\right)\frac{\partial f}{\partial z_{23}}.\\
\end{split}
\end{equation*}

Substituting (\ref{tb}) into the above formula, we conclude Lemma \ref{x12}.\,\,\,$\endpf$
\medskip

Similarly, we can derive that
\begin{lemma}\label{x13}
In the coordinate chart $\left(\left(k_0\cdot\underline{U}_{EI}\right)\times\mathbb R^3,\Lambda_{k_0}\right)$, where $k_0$ is the identity of $SO_3(\mathbb R)$,  \begin{equation*}
\begin{split}
	&-L_{X_{13}}=\left(\frac{1}{2}\sin\alpha\cos\beta\sin 2\gamma+\cos\alpha\tan\beta\sin^2\gamma-\frac{1}{2}\sin\alpha\sin\beta\tan\beta\sin 2\gamma\right)\frac{\partial }{\partial \alpha}\\
		&+\left(\frac{1}{2}\cos\alpha\sin\beta\sin 2\gamma-\sin\alpha\sin^2\beta\cos^2\gamma-\sin\alpha\cos^2\beta\sin^2\gamma\right)\frac{\partial }{\partial \beta}\\
		&+\left(\frac{1}{4}\sin\alpha\sin 2\beta\sin 2\gamma-\frac{1}{2}\sin\alpha\tan\beta\sin 2\gamma+\cos\alpha\sec\beta\sin^2\gamma\right)\frac{\partial }{\partial \gamma}\\
		&+\left(\left(\frac{1}{2}\sin\alpha\sin 2\beta\cos^2\gamma+\frac{1}{2}\sin\alpha\sin 2\beta-\frac{1}{2}\cos\alpha\cos\beta\sin 2\gamma\right)z_{12}\right.\\
		&\,\,\,\,\,\,\,\,\,\,\,\,+\left(\sin\alpha\cos 2\beta\cos\gamma+\cos\alpha\sin\beta\sin\gamma\right)\biggr)\frac{\partial }{\partial z_{12}}\\
		&+\left(\left(\frac{1}{2}\sin\alpha\sin 2\beta \cos 2\gamma-\cos\alpha\cos\beta\sin 2\gamma\right)z_{13}+\left(\sin\alpha\cos 2\beta\cos\gamma\right.\right.\\
		&\,\,\,\,\,\,\,\,\,\,\,\,\left.\left.+\cos\alpha\sin\beta\sin\gamma\right)z_{23}+\left(\frac{1}{2}\sin\alpha \sin 2\beta\sin 2\gamma+\cos\alpha\cos\beta\cos 2\gamma\right)\right)\frac{\partial }{\partial z_{13}}\,\,\,\,\,\,\,\,\,\,\,\,\,\,\,\,\,\,\,\,\,\,\\
		&+\left(\left(-\frac{1}{2}\sin\alpha\sin 2\beta-\frac{1}{2}\sin\alpha\sin 2\beta\sin^2\gamma-\frac{1}{2}\cos\alpha\cos\beta\sin 2\gamma\right)z_{23}\right.\\
		&\,\,\,\,\,\,\,\,\,\,\,\,+\left(\sin\alpha\cos 2\beta\sin\gamma-\cos\alpha\sin\beta\cos\gamma\right)\biggr)\frac{\partial }{\partial z_{23}}.\\
	\end{split}
\end{equation*}
\end{lemma}

\begin{lemma}\label{x23}
In the coordinate chart $\left(\left(k_0\cdot\underline{U}_{EI}\right)\times\mathbb R^3,\Lambda_{k_0}\right)$, where $k_0$ is the identity of $SO_3(\mathbb R)$, 
\small\begin{equation*}
\begin{split}
&-L_{X_{23}}=\left(\cos2\alpha\cos^2\gamma-\cos^2\alpha+\frac{1}{2}\sin 2\alpha\sin\beta\sin 2\gamma\right)\frac{\partial }{\partial \alpha}\\
&+\left(\frac{1}{4}\sin 2\alpha\sin 2\beta\cos 2\gamma-\frac{1}{2}\cos 2\alpha\cos \beta\sin 2\gamma\right)\frac{\partial }{\partial \beta}+\left(-\frac{1}{4}\sin 2\alpha\cos^2\beta\sin 2\gamma\right)\frac{\partial }{\partial \gamma}\\
&+\left(\left(\frac{1}{2}\sin 2\alpha\sin^2\beta \cos^2\gamma-\frac{1}{2}\sin2\alpha\sin^2\gamma-\frac{1}{2}\sin 2\alpha\cos^2\beta-\frac{1}{2}\cos 2\alpha\sin\beta\sin 2\gamma\right)z_{12}\right.\,\,\\
&\,\,\,\,\,\,\,\,\,\,\,\,\left.+\left(\frac{1}{2}\sin 2\alpha\sin 2\beta \cos\gamma-\cos 2\alpha \cos\beta\sin\gamma\right)\right)\frac{\partial }{\partial z_{12}}\\
&+\left(\left(\frac{1}{2}\sin 2\alpha\cos 2\gamma+\frac{1}{2}\sin 2\alpha\sin^2\beta \cos 2\gamma-\cos 2\alpha\sin\beta\sin 2\gamma\right)z_{13}\right.\\
&\,\,\,\,\,\,\,\,\,\,\,\,+\left(\frac{1}{2}\sin 2\alpha\sin 2\beta \cos\gamma-\cos 2\alpha \cos\beta\sin\gamma\right)z_{23}\\
&\,\,\,\,\,\,\,\,\,\,\,\,\left.+\left(\frac{1}{2}\sin 2\alpha \sin 2\gamma+\frac{1}{2}\sin 2\alpha\sin^2\beta\sin 2\gamma+\cos 2\alpha\sin\beta\cos 2\gamma\right)\right)\frac{\partial }{\partial z_{13}}\\
&+\left(\left(\frac{1}{2}\sin 2\alpha\cos^2\beta +\frac{1}{2}\sin 2\alpha\cos^2\gamma-\frac{1}{2}\sin 2\alpha\sin^2\beta \sin^2\gamma-\frac{1}{2}\cos 2\alpha\sin\beta\sin 2\gamma\right)z_{23}\right.\\
&\,\,\,\,\,\,\,\,\,\,\,\,\left.+\left(\frac{1}{2}\sin 2\alpha\sin 2\beta \sin\gamma+\cos 2\alpha\cos\beta\cos\gamma\right)\right)\frac{\partial }{\partial z_{23}}.\\
\end{split}
\end{equation*}
\normalsize
\end{lemma}

\begin{lemma}\label{x21}
In the coordinate chart $\left(\left(k_0\cdot\underline{U}_{EI}\right)\times\mathbb R^3,\Lambda_{k_0}\right)$, where $k_0$ is the identity of $SO_3(\mathbb R)$, 
\begin{equation*}
\begin{split}
&L_{X_{12}}-L_{X_{21}}=\left(\sin\alpha\tan\beta\right)\frac{\partial }{\partial \alpha}+\left(\cos\alpha\right)\frac{\partial }{\partial \beta}+\left(\sin \alpha\sec\beta\right)\frac{\partial }{\partial \gamma}.\\	\end{split}
\end{equation*}
\end{lemma}

\begin{lemma}\label{x31}
In the coordinate chart $\left(\left(k_0\cdot\underline{U}_{EI}\right)\times\mathbb R^3,\Lambda_{k_0}\right)$, where $k_0$ is the identity of $SO_3(\mathbb R)$, 
\begin{equation*}
\begin{split}
&L_{X_{13}}-L_{X_{31}}=-\left(\cos\alpha\tan\beta\right)\frac{\partial }{\partial \alpha}+\left(\sin \alpha\right)\frac{\partial }{\partial \beta}-\left(\cos\alpha\sec\beta\right)\frac{\partial }{\partial \gamma}.\\
\end{split}
\end{equation*}
\end{lemma}

\begin{lemma}\label{x32}
In the coordinate chart $\left(\left(k_0\cdot\underline{U}_{EI}\right)\times\mathbb R^3,\Lambda_{k_0}\right)$, where $k_0$ is the identity of $SO_3(\mathbb R)$, 
\begin{equation*}
\begin{split}
&L_{X_{23}}-L_{X_{32}}=\frac{\partial }{\partial \alpha}.\\
\end{split}
\end{equation*}
\end{lemma}

\subsection{Explicit formulas for the generators of the left-invariant differentials on \texorpdfstring{$SL_3(\mathbb R)$}{dd}}\label{lefti}
In this subsection, we compute the left-invariant differential operators $R\left(E_{ij}\right)$, $1\leq i\neq j\leq 3$ on $SL_3(\mathbb R)$, where $E_{ij}$ is the $3\times 3$ matrix unit with a $1$ in the $i^{\rm th}$ row and $j^{\rm th}$ column. 

Recall the  Euler-Iwasawa coordinates $(\alpha,\beta,\gamma,z_{12}$, $z_{13},z_{23},\lambda_1,\lambda_2)$ of $SL_3(\mathbb R)$ attached to $k_0$ (see Definition \ref{eisl}). Then,
\begin{lemma}\label{r12}
In each coordinate chart $\left(\left(k_0\cdot\underline{U}_{EI}\right)\times\mathbb R^3\times \mathbb R^2,\widehat\Lambda_{k_0}\right)$, 
\begin{equation*}
\begin{split}
&R\left(E_{12}\right)(\alpha,\beta,\gamma,z_{12},z_{13},z_{23},\lambda_1,\lambda_2)=\frac{\lambda_1}{\lambda_2}\frac{\partial}{\partial z_{12}}.\\
\end{split}
\end{equation*}
\end{lemma}
{\noindent\bf Proof of Lemma \ref{r12}.} For matrix-valued functions $g:=\sum_{1\leq i,j\leq 3}g_{ij}(t)E_{ij}$, write $g=O(2)$ if $\left.\frac{d g_{ij}}{dt}(t)\right|_{t=0}=0$, $1\leq i,j\leq 3$.
Computation yields that for each $k_0\in SO_3(\mathbb R)$,
\begin{equation*}
\begin{split}
&k_0\cdot\left(\begin{matrix}k_{11}&k_{12}&k_{13}\\		k_{21}&k_{22}&k_{23}\\		k_{31}&k_{32}&k_{33}\end{matrix}\right)\left(\begin{matrix}1&z_{12}&z_{13}\\		0&1&z_{23}\\
	0&0&1\\\end{matrix}\right)\left(\begin{matrix}\lambda_1&0&0\\	0&\lambda_{2}&0\\
0&0&\lambda_{1}^{-1}\lambda_2^{-1}\end{matrix}\right)\left(\begin{matrix}1&t&0\\		0&1&0\\	0&0&1\\\end{matrix}\right)\\
&=k_0\cdot\left(\begin{matrix} k_{11}& k_{12}& k_{13}\\ k_{21}&k_{22}&k_{23}\\		k_{31}&k_{32}&k_{33}\end{matrix}\right)\left(\begin{matrix}1& z_{12}+\frac{\lambda_1}{\lambda_2}t&z_{13}\\		0&1&z_{23}\\
	0&0&1\\\end{matrix}\right)\left(\begin{matrix}\lambda_1&0&0\\	0&\lambda_{2}&0\\
0&0&\lambda_{1}^{-1}\lambda_2^{-1}\end{matrix}\right)+O(2).\\
\end{split}
\end{equation*}
For smooth functions $f$ on $SL_3(\mathbb R)$, we can derive that
\begin{equation*}
R\left(E_{12}\right)f=\left.\frac{d}{dt}\left(f\left(\alpha,\beta,\gamma,z_{12}+\frac{\lambda_1}{\lambda_2}t,z_{13},z_{23},\lambda_1,\lambda_2\right)\right)\right|_{t=0}=\frac{\lambda_1}{\lambda_2}\frac{\partial f}{\partial z_{12}}.
\end{equation*}

We complete the proof of Lemma \ref{r12}. $\endpf$
\medskip

Similarly, we can derive that
\begin{lemma}\label{r13}
In each coordinate chart $\left(\left(k_0\cdot\underline{U}_{EI}\right)\times\mathbb R^3\times \mathbb R^2,\widehat\Lambda_{k_0}\right)$, 
\begin{equation*}
\begin{split}
&R\left(E_{13}\right)(\alpha,\beta,\gamma,z_{12},z_{13},z_{23},\lambda_1,\lambda_2)=\lambda_1^2\lambda_2\frac{\partial}{\partial z_{13}}.\\
\end{split}
\end{equation*}
\end{lemma}

\begin{lemma}\label{r23} 
In each coordinate chart $\left(\left(k_0\cdot\underline{U}_{EI}\right)\times\mathbb R^3\times \mathbb R^2,\widehat\Lambda_{k_0}\right)$, 
\begin{equation*}
\begin{split}
&R\left(E_{23}\right)(\alpha,\beta,\gamma,z_{12},z_{13},z_{23},\lambda_1,\lambda_2)=\lambda_1\lambda_2^2\left(\frac{\partial}{\partial z_{23}}+z_{12}\frac{\partial}{\partial z_{13}}\right).\\
\end{split}
\end{equation*}
\end{lemma}
\begin{lemma}\label{r21} 
In each coordinate chart $\left(\left(k_0\cdot\underline{U}_{EI}\right)\times\mathbb R^3\times \mathbb R^2,\widehat\Lambda_{k_0}\right)$, 
\begin{equation*}
\begin{split}
&R\left(E_{21}\right)(\alpha,\beta,\gamma,z_{12},z_{13},z_{23},\lambda_1,\lambda_2)=\frac{\lambda_2}{\lambda_1}\left\{\sec\beta\sin\gamma\frac{\partial }{\partial \alpha}+\cos\gamma\frac{\partial }{\partial \beta}+\tan\beta\sin\gamma\frac{\partial}{\partial \gamma}\right.\\
&\,\,\,\,\,\,\,\,\,\,\left.+\left(z_{12}^2+1\right)\frac{\partial}{\partial z_{12}}+z_{23}\frac{\partial}{\partial z_{13}}-z_{13}\frac{\partial}{\partial z_{23}}\right\}+\lambda_2z_{12}\frac{\partial}{\partial \lambda_1}-\frac{\lambda_2^2}{\lambda_1}z_{12}\frac{\partial}{\partial \lambda_2}.\\
\end{split}
\end{equation*}
\end{lemma}

\begin{lemma}\label{r31} 
In each coordinate chart $\left(\left(k_0\cdot\underline{U}_{EI}\right)\times\mathbb R^3\times \mathbb R^2,\widehat\Lambda_{k_0}\right)$, 
\begin{equation*}
\begin{split}
&R\left(E_{31}\right)(\alpha,\beta,\gamma,z_{12},z_{13},z_{23},\lambda_1,\lambda_2)=\lambda_1^{-2}\lambda_2^{-1}\left\{\left(\left(\sec\beta\sin\gamma\right) z_{23}-\left(\sec\beta\cos\gamma\right) z_{12}\right)\frac{\partial}{\partial \alpha}\right.\\
&\,\,\,\,\,\,\,\,\,\,+\left(\left(\cos\gamma\right) z_{23}+\left(\sin\gamma\right) z_{12}\right)\frac{\partial}{\partial \beta}+\left(\left(\tan\beta\sin\gamma\right) z_{23}-\left(\tan\beta\cos\gamma\right) z_{12}-1\right) \frac{\partial}{\partial\gamma}\\
&\,\,\,\,\,\,\,\,\,\,\left.+\left(z_{23}+z_{12}^2z_{23}\right)\frac{\partial}{\partial z_{12}}+\left(1+z_{13}^2+z_{23}^2-z_{12}z_{13}z_{23}\right)\frac{\partial}{\partial z_{13}}+\left(-z_{12}-z_{12}z_{23}^2\right)\frac{\partial}{\partial z_{23}}\right\}\\
&\,\,\,\,\,\,\,\,\,\,+\lambda_1^{-1}\lambda_2^{-1}z_{13}\frac{\partial}{\partial \lambda_1}-\lambda_1^{-2}z_{12}z_{23}\frac{\partial}{\partial \lambda_2}.\\
\end{split}
\end{equation*}
\end{lemma}

\begin{lemma}\label{r32} 
In each coordinate chart $\left(\left(k_0\cdot\underline{U}_{EI}\right)\times\mathbb R^3\times \mathbb R^2,\widehat\Lambda_{k_0}\right)$,
\begin{equation*}
\begin{split}
&R\left(E_{32}\right)(\alpha,\beta,\gamma,z_{12},z_{13},z_{23},\lambda_1,\lambda_2)=\lambda_1^{-1}\lambda_2^{-2}\left\{\sec\beta\cos\gamma\frac{\partial }{\partial \alpha}-\sin\gamma\frac{\partial }{\partial \beta}+\tan\beta\cos\gamma\frac{\partial}{\partial \gamma}\right.\\
&\,\,\,\,\,\,\,\,\,\left.+\left(z_{13}-z_{12}z_{23}\right)\frac{\partial}{\partial z_{12}}+z_{13}z_{23}\frac{\partial}{\partial z_{13}}+\left(z_{23}^2+1\right)\frac{\partial}{\partial z_{23}}\right\}+\lambda_1^{-1}\lambda_2^{-1}z_{23}\frac{\partial}{\partial \lambda_2}.\\	\end{split}
\end{equation*}
\end{lemma}

\section{Regularization}\label{REG}

{\noindent\bf Proof of Lemma \ref{moll}.}
For each $u\in U(\mathfrak{sl}_3(\mathbb{R}))$, it holds by definition that $L_u(\psi_m*f)=L_u(\psi_m)*f$. This implies (a) since  $\{L_u\}_{u\in U(\mathfrak{sl}_3(\mathbb R))}$ provide all partial derivatives with respect to any local coordinates at each point of the manifold $X$. 
By the Cauchy-Schwarz inequality, 
\begin{equation*}
\begin{split}
&\int_{X}|(\psi_m*f)(x)|^2{\rm d}\mu=\int_{X}\left|\int_{SL_3(\mathbb R)}\psi_m(g)\cdot f\left(g^{-1}x\right){\rm d}g\right|^2{\rm d}\mu\\  
&\leq\int\limits_{X}\left(\int\limits_{SL_3(\mathbb R)}\psi_m\left(\widetilde g\right) {\rm d}\widetilde g\right) \left(\int\limits_{SL_3(\mathbb R)}\psi_m(g)\cdot \left|f\left(g^{-1}x\right)\right|^2{\rm d}g\right){\rm d}\mu= \int_{X} \left|f\left(x\right)\right|^2{\rm d}\mu.\\
\end{split} 
\end{equation*}
For each $\epsilon>0$, we can find continuous function $f_\epsilon$ with compact support such that 
\begin{equation}\label{f1}
\int_{X} \left|f_{\epsilon}\left(x\right)-f(x)\right|^2{\rm d}\mu<\epsilon.  \end{equation}
By the Cauchy-Schwarz inequality again, 
\begin{equation*}
\int_{X}|(\psi_m*f_{\epsilon})(x)-f_{\epsilon}(x)|^2{\rm d}\mu\leq\sup_{g\in U_m}\int_{X} \left|f_{\epsilon}\left(g^{-1}x\right)-f_{\epsilon}(x)\right|^2{\rm d}\mu. 
\end{equation*}
Since $f_{\epsilon}$ is uniform continuous with compact support, we conclude \begin{equation}\label{152}
\lim_{m\rightarrow\infty}\int_{X}|(\psi_m*f_{\epsilon})(x)-f_{\epsilon}(x)|^2{\rm d}\mu=0. \end{equation} 
Then, we can derive that
\begin{equation*}
\lim_{m\rightarrow\infty}\int_{X}|(\psi_m*f)(x)-f(x)|^2{\rm d}\mu=0.
\end{equation*}
(b) is now proved. The same argument yields (c) for $L_u(\psi_m*f)=(L_u\psi_m)*f$.

Next, assume $f,\Delta f\in L^2(X)$. We first show that $\Delta(\psi_m*f)=\psi_m*(\Delta f)$ in the sense of distribution. Recall the notion that $\{X_{ij}\}_{1\leq i,j\leq 3,\,(i,j)\neq (3,3)}$is a basis of $\mathfrak{sl}_3(\mathbb R)$ (see (\ref{slba})) and that the lift $\widetilde w$ of $w\in C^{\infty}\left(SL_3(\mathbb R)/A\right)$ is defined by $\widetilde w:=w\circ\pi$ (see \eqref{hgh}). For clarity, in the following, we denote the points of $X$ by $xA$. 

For each invariant differential $\Delta$, there is a polynomial $P_{\Delta}$ by Theorem \ref{HCiso} such that
\begin{equation*}
\begin{split}
&(\Delta w)(xA)=\left.P_{\Delta}\left(\cdots,\frac{\partial}{\partial t_{ij}},\cdots\right)\widetilde w\left(x\exp\left(\sum_{\substack{1\leq i,j\leq 3\\(i,j)\neq (3,3)}}t_{ij}X_{ij}\right)\right)\right|_{t_{11}=t_{12}=\cdots=t_{32}=0}\\
&\,\,\,\,\,\,\,\,\,\,\,\,\,\,\,\,\,\,\,\,\,=:\left.P_{\Delta}\left(\frac{\partial}{\partial t_{ij}}\right)\widetilde w\left(x\exp\left(\sum t_{ij}X_{ij}\right)\right)\right|_{t_{ij}=0},\,\,\,\,\,\,\,\,w\in C^{\infty}\left(SL_3(\mathbb R)/A\right).
\end{split}
\end{equation*}
Write $\psi_m^*(g):=\psi_m(g^{-1})$. Noticing that the formal adjoint $\Delta^*$ of $\Delta$ is invariant on $X$ (\cite{Ba}), for each $h\in C^{\infty}_c(X)$, we can derive that
\begin{equation*}
\begin{split}
&(\Delta(\psi_m*f),h)=\int_{X}f\left(xA\right)\overline{P_{\Delta^*  }\left(\frac{\partial}{\partial t_{ij}}\right)\left.\left( \int_{SL_3(\mathbb R)}\psi_m^*(g)\widetilde h\left(g^{-1}x\exp\left(\sum t_{ij}X_{ij}\right)\right)\,{\rm d}g\right)\right|_{t_{ij}=0}},\\
&=\int_{X}(\Delta f)(xA)\cdot\left(\int_{SL_3(\mathbb R)}\psi_m\left(g^{-1}\right)\cdot\overline{h\left(g^{-1}xA\right)}\,{\rm d}g\right)\,{\rm d}\mu\\
&=\int_{X}\left(\int_{SL_3(\mathbb R)}\psi_m\left(\widetilde g\right)\cdot(\Delta f)\left(\widetilde g^{-1}\widetilde xA\right)\,{\rm d}\widetilde g\right)\cdot\overline{h\left(\widetilde xA\right)}\,{\rm d}\mu=(\psi_m*(\Delta f),h).\\
\end{split}
\end{equation*}
Then,
\begin{equation}
\begin{split}
&|(\Delta(\psi_m*f),h)|=|(\psi_m*(\Delta f),h)|\leq\int_{X}\left|(\Delta f)(x)\right|^2\,{\rm d}\mu\cdot\int_{X}\left|h(x)\right|^2\,{\rm d}\mu.\\    
\end{split}
\end{equation}
Now $\psi_m*\Delta f$ lies in $L^2(X)$ thanks to (b), we also have $\Delta(\psi_m*f)=\psi_m*\Delta f$ in $L^2(X)$. Using (b) again, 
$\Delta(\psi_m*f)=\psi_m*\Delta f\to \Delta f$ in $L^2(X)$ as $m\to\infty$. 
We complete the proof of Lemma \ref{moll}.\,\,\,$\endpf$

\end{appendices}

\bibliographystyle{plain}
\bibliography{sn-bibliography}%

\begin{thebibliography}{10}

\bibitem{BK}
Yves Benoist and Toshiyuki Kobayashi.
\newblock Tempered reductive homogeneous spaces.
\newblock {\em J. Eur. Math. Soc. (JEMS)}, 17(12):3015--3036, 2015.

\bibitem{BKS}
Alexey Bondal, Mikhail Kapranov, and Vadim Schechtman.
\newblock Perverse schobers and birational geometry.
\newblock {\em Selecta Math. (N.S.)}, 24(1):85--143, 2018.

\bibitem{Ca}
Stefan Catoiu.
\newblock Prime ideals of the enveloping algebra {$U({\rm sl}_3)$}.
\newblock {\em Comm. Algebra}, 28(2):981--1027, 2000.

\bibitem{Ch}
Paul~R. Chernoff.
\newblock Essential self-adjointness of powers of generators of hyperbolic equations.
\newblock {\em J. Functional Analysis}, 12:401--414, 1973.

\bibitem{Co}
H.~O. Cordes.
\newblock Self-adjointness of powers of elliptic operators on non-compact manifolds.
\newblock {\em Math. Ann.}, 195:257--272, 1972.

\bibitem{De}
Patrick Delorme.
\newblock Formule de {P}lancherel pour les espaces sym\'{e}triques r\'{e}ductifs.
\newblock {\em Ann. of Math. (2)}, 147(2):417--452, 1998.

\bibitem{DKKS}
Patrick Delorme, Friedrich Knop, Bernhard Kr\"{o}tz, and Henrik Schlichtkrull.
\newblock Plancherel theory for real spherical spaces: construction of the {B}ernstein morphisms.
\newblock {\em J. Amer. Math. Soc.}, 34(3):815--908, 2021.

\bibitem{DZ}
Anthony~H. Dooley and Genkai Zhang.
\newblock Spherical functions on harmonic extensions of {$H$}-type groups.
\newblock {\em J. Geom. Anal.}, 9(2):247--255, 1999.

\bibitem{FJ}
Mogens Flensted-Jensen.
\newblock Discrete series for semisimple symmetric spaces.
\newblock {\em Ann. of Math. (2)}, 111(2):253--311, 1980.

\bibitem{Ga}
Matthew~P. Gaffney.
\newblock A special {S}tokes's theorem for complete {R}iemannian manifolds.
\newblock {\em Ann. of Math. (2)}, 60:140--145, 1954.

\bibitem{HC1}
Harish-Chandra.
\newblock Spherical functions on a semisimple {L}ie group. {I}.
\newblock {\em Amer. J. Math.}, 80:241--310, 1958.

\bibitem{HC2}
Harish-Chandra.
\newblock Spherical functions on a semisimple {L}ie group. {II}.
\newblock {\em Amer. J. Math.}, 80:553--613, 1958.

\bibitem{He}
Sigurdur Helgason.
\newblock {\em Groups and geometric analysis}, volume~83 of {\em Mathematical Surveys and Monographs}.
\newblock American Mathematical Society, Providence, RI, 2000.
\newblock Integral geometry, invariant differential operators, and spherical functions, Corrected reprint of the 1984 original.

\bibitem{HY}
Xiaojun Huang and Wanke Yin.
\newblock Flattening of {CR} singular points and analyticity of the local hull of holomorphy {II}.
\newblock {\em Adv. Math.}, 308:1009--1073, 2017.

\bibitem{KK}
Fanny Kassel and Toshiyuki Kobayashi.
\newblock Poincar\'{e} series for non-{R}iemannian locally symmetric spaces.
\newblock {\em Adv. Math.}, 287:123--236, 2016.

\bibitem{Kn}
Friedrich Knop.
\newblock A {H}arish-{C}handra homomorphism for reductive group actions.
\newblock {\em Ann. of Math. (2)}, 140(2):253--288, 1994.

\bibitem{MOZ}
J.~M\"{o}llers, B.~\O~rsted, and G.~Zhang.
\newblock Invariant differential operators on {H}-type groups and discrete components in restrictions of complementary series of rank one semisimple groups.
\newblock {\em J. Geom. Anal.}, 26(1):118--142, 2016.

\bibitem{NS}
Edward Nelson and W.~Forrest Stinespring.
\newblock Representation of elliptic operators in an enveloping algebra.
\newblock {\em Amer. J. Math.}, 81:547--560, 1959.

\bibitem{OZ}
Bent \O~rsted and Gen~Kai Zhang.
\newblock Weyl quantization and tensor products of {F}ock and {B}ergman spaces.
\newblock {\em Indiana Univ. Math. J.}, 43(2):551--583, 1994.

\bibitem{OS}
Toshio Oshima and Jiro Sekiguchi.
\newblock Eigenspaces of invariant differential operators on an affine symmetric space.
\newblock {\em Invent. Math.}, 57(1):1--81, 1980.

\bibitem{PZ}
Lizhong Peng and Genkai Zhang.
\newblock Tensor products of holomorphic representations and bilinear differential operators.
\newblock {\em J. Funct. Anal.}, 210(1):171--192, 2004.

\bibitem{Rob}
Joel Roberts.
\newblock Old and new results about the triangle varieties.
\newblock In {\em Algebraic geometry ({S}undance, {UT}, 1986)}, volume 1311 of {\em Lecture Notes in Math.}, pages 197--219. Springer, Berlin, 1988.

\bibitem{RSe}
Joel Roberts and Robert Speiser.
\newblock Enumerative geometry of triangles. {I}.
\newblock {\em Comm. Algebra}, 12(9-10):1213--1255, 1984.

\bibitem{Ros}
Wulf Rossmann.
\newblock Analysis on real hyperbolic spaces.
\newblock {\em J. Functional Analysis}, 30(3):448--477, 1978.

\bibitem{SV}
Yiannis Sakellaridis and Akshay Venkatesh.
\newblock Periods and harmonic analysis on spherical varieties.
\newblock {\em Ast\'{e}risque}, (396):viii+360, 2017.

\bibitem{Sc}
H.~Schubert.
\newblock Anzahlgeometrische {B}ehandlung des {D}reiecks.
\newblock {\em Math. Ann.}, 17(2):153--212, 1880.

\bibitem{Seg}
I.~E. Segal.
\newblock Hypermaximality of certain operators on {L}ie groups.
\newblock {\em Proc. Amer. Math. Soc.}, 3:13--15, 1952.

\bibitem{Se}
J.~G. Semple.
\newblock The triangle as a geometric variable.
\newblock {\em Mathematika}, 1:80--88, 1954.

\bibitem{Sh1}
Goro Shimura.
\newblock Arithmetic of differential operators on symmetric domains.
\newblock {\em Duke Math. J.}, 48(4):813--843, 1981.

\bibitem{Sh2}
Goro Shimura.
\newblock Differential operators and the singular values of {E}isenstein series.
\newblock {\em Duke Math. J.}, 51(2):261--329, 1984.

\bibitem{Sh3}
Goro Shimura.
\newblock On differential operators attached to certain representations of classical groups.
\newblock {\em Invent. Math.}, 77(3):463--488, 1984.

\bibitem{Sh4}
Goro Shimura.
\newblock Invariant differential operators on {H}ermitian symmetric spaces.
\newblock {\em Ann. of Math. (2)}, 132(2):237--272, 1990.

\bibitem{Sh5}
Goro Shimura.
\newblock Differential operators, holomorphic projection, and singular forms.
\newblock {\em Duke Math. J.}, 76(1):141--173, 1994.

\bibitem{St}
Robert~S. Strichartz.
\newblock Harmonic analysis on hyperboloids.
\newblock {\em J. Functional Analysis}, 12:341--383, 1973.

\bibitem{BS1}
E.~P. van~den Ban and H.~Schlichtkrull.
\newblock The {P}lancherel decomposition for a reductive symmetric space. {I}. {S}pherical functions.
\newblock {\em Invent. Math.}, 161(3):453--566, 2005.

\bibitem{BS2}
E.~P. van~den Ban and H.~Schlichtkrull.
\newblock The {P}lancherel decomposition for a reductive symmetric space. {II}. {R}epresentation theory.
\newblock {\em Invent. Math.}, 161(3):567--628, 2005.

\bibitem{Ba}
Erik~P. van~den Ban.
\newblock Invariant differential operators on a semisimple symmetric space and finite multiplicities in a {P}lancherel formula.
\newblock {\em Ark. Mat.}, 25(2):175--187, 1987.

\bibitem{W}
Joseph~A. Wolf.
\newblock Essential self-adjointness for the {D}irac operator and its square.
\newblock {\em Indiana Univ. Math. J.}, 22:611--640, 1972/73.

\bibitem{ZG3}
Genkai Zhang.
\newblock Invariant differential operators on {H}ermitian symmetric spaces and their eigenvalues.
\newblock {\em Israel J. Math.}, 119:157--185, 2000.

\bibitem{ZG2}
Genkai Zhang.
\newblock Shimura invariant differential operators and their eigenvalues.
\newblock {\em Math. Ann.}, 319(2):235--265, 2001.

\bibitem{ZG1}
GenKai Zhang.
\newblock Tensor products of complementary series of rank one {L}ie groups.
\newblock {\em Sci. China Math.}, 60(11):2337--2348, 2017.

\bibitem{ZZ1}
Han Zhang and Runlin Zhang.
\newblock Nondivergence on homogeneous spaces and rigid totally geodesic submanifolds.
\newblock {\em Israel J. Math.}, 264(1):149--176, 2024.

\bibitem{Zh1}
Runlin Zhang.
\newblock Limiting distribution of translates of the orbit of a maximal {$\mathbb Q$}-torus from identity on {${\rm SL}_N(\mathbb R)/{\rm SL}_N(\mathbb Z)$}.
\newblock {\em Math. Ann.}, 375(3-4):1231--1281, 2019.

\bibitem{Zh2}
Runlin Zhang.
\newblock Translates of homogeneous measures associated with observable subgroups on some homogeneous spaces.
\newblock {\em Compos. Math.}, 157(12):2657--2698, 2021.

\bibitem{Zh3}
Runlin Zhang.
\newblock Equidistribution of translates of a homogeneous measure on the {B}orel-{S}erre compactification.
\newblock {\em Discrete Contin. Dyn. Syst.}, 42(4):2053--2071, 2022.

\bibitem{Zh4}
Runlin Zhang.
\newblock Asymptotics of integral points, equivariant compactifications and equidistributions for homogeneous spaces.
\newblock {\em arXiv preprint arXiv:2408.02325}, 2024.

\end{thebibliography}

\end{document}